\documentclass{article}
\usepackage{amsmath,amsfonts,amsthm,amscd,amssymb,latexsym,comment}
\usepackage{enumerate,graphicx,psfrag,subfigure,changebar,oldgerm}
\title{Centraliser Dimension of Partially Commutative Groups
\footnote{Research supported by EPSRC grant GR/R87406/01 }
\footnote{{\em{2000 Mathematics Subject Classification}} 20F36 (primary), 20E15, 20F10 (secondary)}
}
\author{ \textsf{Andrew J. Duncan}
 \and \textsf{Ilya V. Kazachkov\footnote{Supported by RFFI grant N05-01-00057-a}}
 \and \textsf{Vladimir N. Remeslennikov\footnotemark[3]}
}
\def\nul{\emptyset }

\def\b{\beta }

\def\e{\varepsilon }
\def\G{\Gamma }

\def\a{\alpha }

\def\fC{{\textswab C}}
\newtheorem{thm}{Theorem}[section]
\newtheorem{lem}[thm]{Lemma}
\newtheorem{cor}[thm]{Corollary}

\newtheorem{defn}[thm]{Definition}
\newtheorem{exam}[thm]{Example}
\newenvironment{expl}{\begin{exam} \rm}{\end{exam}}
\numberwithin{equation}{section}

\newcommand{\cd}{\texttt{cdim}}

\newcommand{\az}{\mathop{{\ensuremath{\alpha}}}}
\newcommand{\la}{\langle}
\newcommand{\ra}{\rangle}
\newcommand{\sdc}{>}
\newcommand{\sac}{<}

\newcommand{\edc}{\ge}
\newcommand{\eac}{\le}

\newcommand{\neac}{\nleq}

\newcommand{\be}{\begin{enumerate}}
\newcommand{\ee}{\end{enumerate}}

\begin{document}
\maketitle
\begin{abstract}
In paper \cite{DKR} we investigated the centraliser dimension of
groups. In the current paper we study properties of centraliser
dimension for the class of free partially commutative groups and, as a
corollary, we obtain an efficient algorithm for computation of
centraliser dimension in these groups.
\end{abstract}
\section{Introduction}

Partially commutative groups have been extensively studied
in several different guises. They are variously known as
right-angled Artin groups; trace groups;
graph groups or  locally free groups. For an introduction to this
class of groups and a survey of the literature see \cite{BEKR}.

In this paper we study the lattices of centralisers in partially
commutative groups. 
Centraliser lattices, of  various
groups, have been investigated by numerous authors: for example
\cite{Ito},  
\cite{Zaleski65}, \cite{Vas}, \cite{LR}, \cite{Schmidt94},
\cite{Schmidt95}, \cite{Weh}, 
\cite{Kos73}, \cite{kegel73}, \cite{BryantHartley79}, \cite{Bryant79},
\cite{AS}, \cite{Bl}, \cite{Wag}, \cite{MS} and \cite{DKR}.
In particular, a detailed account of results in the field
can be found in V.~A.~Antonov's book  \cite{book}. 
Following \cite{Schmidt94} and \cite{MS}, here we focus attention 
on the height of the lattice of centralisers of a group: which
we call the centraliser dimension.
Our interest in centralisers stems from their importance in
constructing algebraic geometry over 
groups. We plan to apply the results of this paper to the
study of, in the first instance, solutions of equations in
one variable in partially commutative groups.

Partially commutative groups have been extensively studied
in several different guises. They are variously known as
right-angled Artin groups; trace groups;
graph groups or  locally free groups. For an introduction to this
class of groups and a survey of the literature see \cite{BEKR}.

Let $\G$ be a finite, undirected, simple graph. Let $X=V(\G)$ be the set of vertices
of $\G$ and let $F(X)$ be the free group on $X$.
For elements $g,h$ of a group we denote the commutator $g^{-1}h^{-1}gh$
of $g$ and $h$ by $[g,h]$. Let
\[
R=\{[x_i,x_j]\in F(X)\mid x_i,x_j\in X \textrm{ and there is an edge from } x_i \textrm{ to } x_j
\textrm{ in } \G\}.
\]
We define
the {\em free partially commutative group with (commutation) graph } $\G$ to be the
group $G(\G)$ with presentation
$
\left< X\mid R\right>.
$

The elements of $X$ are termed the {\em canonical} generators of $G(\G)$. In this
article we refer to
finitely generated free partially commutative groups as {\em partially
  commutative} groups.

If $S$ is a subset of a group $G$ then the centraliser
of $S$ in $G$ is
$C_G(S)=\{g\in G: gs=sg, \textrm{ for all } s\in S\}$. We write
$C(S)$ instead of $C_G(S)$ when the meaning is clear.
Let $\fC(G)$ denote the set of centralisers of a group $G$.
The relation of inclusion then defines  a partial order `$\le$' on $\fC(G)$.
We define the infimum of a pair of elements of $\fC(G)$ as
obvious way:
    \[
        C(M_1) \wedge C(M_2)= C(M_1)\cap C(M_2)=C( M_1 \cup M_2).
    \]
Moreover  the supremum $C(M_1) \vee C(M_2)$ of elements
$C(M_1)$ and $C(M_2)$ of $\fC(G)$ may
be
defined to be the intersection of all centralisers containing $C(M_1)$ and $C(M_2)$.
Then  $C(M_1) \vee C(M_2)$ is  minimal among centralisers containing $C(M_1)$ and $C(M_2)$.
These definitions make $\fC(G)$ into a lattice, called the
{\em centraliser lattice} of $G$. This lattice is bounded as it has
a greatest element, $G= C(1)$, and a least element, $Z(G)$, the centre
of $G$. 

If $C$ and
$C^\prime$
are in $\fC(G)$ with $C$ strictly contained in $C^\prime$ we
write $C<C^\prime$.
If $C_i$ is a centraliser, for $i=0,\ldots ,k$, with $C_0>\cdots >C_k$
then we call $C_0,\ldots, C_k$ a {\em centraliser chain} of {\em
  length} $k$. Infinite descending, ascending and doubly-infinite centraliser
chains are defined in the obvious way.
\begin{defn}[{\em cf.} \cite{Schmidt94}, \cite{MS}] \label{def:cd}
If there exists an integer $d$ such that
the group $G$  has a centraliser chain of length $d$ and no centraliser
 chain of length greater than $d$ then $G$ is said to have
{\em centraliser dimension} $\cd (G)=d$.
If no such integer $d$ exists
we define $\cd (G)= \infty$.
\end{defn}
If $\cd(G)=d$ then every strictly descending chain of centralisers in $\fC(G)$
from $G$
to $Z(G)$ contains at most $d$ inclusions. This number
is usually referred to as the {\em height} of the lattice;
so
$\cd (G)$ is the height of the centraliser lattice of $G$.

 Throughout the remainder of this paper
$G$ denotes a (finitely generated free) partially commutative
group (unless we explicitly indicate otherwise).
The centraliser dimension of partially commutative groups
is finite because
all partially commutative groups are linear \cite{H} 
 and all linear groups have finite
centraliser dimension \cite{MS}.

The centraliser dimension of a partially commutative group is easy to calculate
and depends only on the centralisers of subsets of $X$. If $Y\subseteq X$ then
we call $C(Y)$ a {\em canonical} centraliser. From the above it follows that
the intersection of two canonical centralisers is again a canonical centraliser.
As shown in Section \ref{sec:cancent} below the supremum, in $\fC(G)$,
of two canonical centralisers
is also a canonical centraliser. Hence the set $\fC(X;G)$ of canonical centralisers
forms a sublattice of $\fC(G)$. We show in Theorem \ref{thm:pcdimgen} that
the centraliser dimension of $G$ is equal to the height of the lattice  $\fC(X;G)$.

Now let $x$ be a vertex of $\G$ and let $\G_x$ denote the graph formed
from $\G$ by removal of the vertex $x$. Let $G=G(\G)$ and let $G_x=G(\G_x)$, the
partially commutative group with graph $\G_x$. In Theorem \ref{cd-ext} we show that
$\cd(G)=\cd(G_x)+d$, where $d=0,1$ or $2$. Moreover the theorem describes exactly
how $d$ depends on $\G$ and $x$. We conclude the paper by applying
Theorem \ref{cd-ext} to calculate the centraliser
dimension of groups in some specific classes of partially commutative groups.

The authors thank Claas R\"over for helpful conversations and careful reading
of, in particular Theorem \ref{thm:pcdimgen}.

\section{Preliminaries}\label{sec:prelim}

First we recall some of the notation and definitions of \cite{BEKR}.
If $w\in (X\cup X^{-1})^*$ then we denote by $\a(w)$ the set of
elements $x\in X$ such that $x$ or $x^{-1}$ occurs in $w$. By abuse of
notation we identify words of $(X\cup X^{-1})^*$ with the elements of
$G$ which they represent. Moreover for $u,v \in (X\cup X^{-1})^*$ we
use $u=v$ to mean $u$ and $v$ are equal as elements of $G$. Equality
of words $u,v$ in $(X\cup X^{-1})^*$ is denoted by $u\equiv v$.
If $g\in G$ and $w\in (X\cup X^{-1})^*$ is a
word of minimal length representing $g$ then we say that $w$ is a {\em
  minimal} form of $g$ (or just that $w$ is {\em minimal}).

The Transformation Lemma, \cite[Lemma 2.3]{BEKR} asserts that, if $u$
and $v$ are minimal and $u=v$ then we can transform the word $u$ into
the word $v$ using only commutation relations from $R$ (that is,
without insertion or deletion of any subwords of the form
$x^{\e}x^{-\e}, x\in X$).
From the
Transformation Lemma it follows that if $u$
and $v$ are minimal and $u=v$ then $\a(u)=\a(v)$ and $|u|=|v|$
(where $|w|$ denotes the length of the word $w$). Therefore, for
any element $g\in G$ we may define $\a(g)=\a(w)$, and the {\em length} $l(g)$ of $g$ as
$l(g)=|w|$, where $w$ is a
minimal form of $g$.

If $g,h\in G$ such that $l(gh)=l(g)+l(h)$ then we write $gh=g\circ h$.
It follows that $gh=g\circ h$ if and only if, for all minimal forms
$u$ and $v$ of $g$ and $h$, respectively, $uv$ is a minimal form for
$gh$. Clearly if $w\in (X\cup X^{-1})^*$ is minimal and $w\equiv uv$
then, in $G$, $w=u\circ v$. If $k=g\circ h$ then we say that $g$ is a
{\em left divisor} of $k$ (and $h$ is a {\em right divisor} of $k$).
We denote by $D_l(k)$ and $D_r(k)$ the set of left and right divisors of $k$,
respectively.
We shall use the following Lemma ({\em cf.} \cite[Lemma 3.9]{BEKR}).
\begin{lem}\label{lem:fl}
If $a\circ b=c\circ d$. Then $a=c_1\circ d_1$ and $b=c_2\circ d_2$, where
$c=c_1\circ c_2$, $d=d_1\circ d_2$ and $[\a(c_2),\a(d_1)]=1$.
\end{lem}
\begin{proof}
If $l(a)=0$ then the result holds, trivially. Suppose then that $l(a)>0$ and that
the result holds for all $a$ of shorter length. Then we may write
$a=x\circ a^\prime$, where $x\in X^{\pm 1}$. In this case $x\in D_l(c\circ d)$, 
so either $x\in D_l(c)$ or $x\in D_l(d)$ and $[x,\a(c)]=1$. Assume first that 
the latter holds. Then $d=x\circ d^\prime$ and it follows that $a^\prime \circ b
=c\circ d^\prime$. As $l(a^\prime)<l(a)$ there exist $c_1$, $c_2$, $d^\prime_1$
and $d_2$ such that $a^\prime=c_1\circ d^\prime_1$, $b=c_2\circ d_2$,
 $c=c_1\circ c_2$, $d^\prime=d_1^\prime \circ d_2$ and 
$[\a(c_2),\a(d^\prime_1)]=1$. The result now follows on setting $d_1=x\circ d_1^\prime$.
In the case where $x\in D_l(c)$ a similar, simpler, argument shows that the lemma
holds.
\end{proof} 

Let $u$ and $v$ be elements of $G$. A {\em greatest common left divisor} of $u$ and
$v$ is an element $c\in D_l(u)\cap D_l(v)$ such that every element of
$D_l(u)\cap D_l(v)$ is a left divisor of $c$. In \cite{BEKR} it is shown that
every pair of elements $u$, $v$ of $G$ has a unique greatest common left
divisor, which we denote by $\gcd_l(u,v)$. This definition extends to $n$-tuples
of elements of $G$; we recursively define $\gcd_l(u_1,\ldots,u_n)$ to be
$\gcd_l(\gcd_l(u_1,\ldots, u_{n-1}),u_n)$. It follows from \cite{BEKR} that
$\gcd_l(u_1,\ldots,u_n)$  is independent of the order of $u_1,\ldots ,u_n$.
The {\em greatest common right divisor}
$\gcd_r(u_1,\ldots,u_n)$
of an $n$-tuple of elements is defined analogously.

A subgroup of $G$ generated by a subset $Y$ of $G$ is called an {\em
canonical parabolic} subgroup. 
Let $Y\subseteq X$ and denote by $G(Y)$ the
partially commutative group given by the presentation with
generators $Y$ and relators those $[x_i,x_j]\in R$ such that
$x_i,x_j\in Y$. Thus $G(Y)$ has commutation graph the full subgraph
of $\G(G)$ generated by $Y$. It follows from the Transformation
Lemma, \cite[Lemma 2.5]{BEKR}, that $G(Y)$ is the canonical
parabolic subgroup of $G$ generated by $Y$. If $g\in G$ and $H$ is
a canonical parabolic subgroup of $G$ then we say that $H^g$ is a
{\em parabolic subgroup}. Thus parabolic subgroups are those with
generating set $\{y^g:g\in Y\}$, for some subset $Y$ of $X$ and some
$g\in G$.

As in \cite{BEKR}, given an element $g\in G$ and a subset $Y$ of $X$,
the set $D_l(g,Y)$ of
{
\em parabolic left divisors
}
of
$g$ in $G(Y)$ is defined to be $D_l(g)\cap G(Y)$. A
{\em greatest common left divisor}
$\gcd_l(g,Y)$
 of $g$ in $G(Y)$ is an element $d$
of $D_l(g,Y)$ such that every element of $D_l(g,Y)$
is a left divisor of $d$. It follows from \cite{BEKR} that
$\gcd_l(g,Y)$ exists and is unique,
for
all elements $g\in G$ and all canonical parabolic subgroups $G(Y)$ of $G$.

We say that $w\in (X\cup X^{-1})^*$ is {\em cyclically minimal}
if $l(w)\le l(w^g)$, for all $g\in G$.
If $w$ is minimal and $w\equiv uv$ with $u$ and $v$ minimal then we
call  $vu\in (X\cup X^{-1})^*$ a {\em cyclic permutation} of $w$.
\begin{lem}\label{cycgeod}
The following are equivalent.
\be[(i)]
\item\label{cycgeod1} $w$ is cyclically minimal.
\item\label{cycgeod2} If $y\in X\cup X^{-1}$ is a left divisor of $w$ then $y^{-1}$ is not a right
  divisor of $w$.
\item\label{cycgeod3} All cyclic permutations of $w$ are minimal.
\ee
Moreover, if $w$ is cyclically minimal then $\a(w^g)\supseteq \a(w)$,
for all $g\in G$.
\end{lem}
\begin{proof}
The equivalence of (\ref{cycgeod1}) and (\ref{cycgeod2}) is
Proposition 5.6 of \cite{BEKR}. To see that (\ref{cycgeod1}) implies
(\ref{cycgeod3}) suppose that $w$ is cyclically minimal and that $vu$
is a cyclic permutation of $w$. Then $vu=vwv^{-1}$ so $l(vu)\ge
l(w)=l(u)+l(v)$. Hence $vu$ is minimal as required. To  complete the
proof of the equivalence we show that (\ref{cycgeod3}) implies (\ref{cycgeod2}). Suppose
that $y\in X^{\pm 1}$ is a left divisor of $w$ and that $y^{-1}$ is a
right divisor of $w$. Then $w=y\circ w_1\circ y^{-1}$, so the cyclic
permutation $w_1 y^{-1}y$ is not minimal, a contradiction.

For the final statement suppose that $w$ is cyclically minimal and
that $u$ is a minimal form for an element of $G$. The result holds
trivially if $u=1$ and we use induction on the length of $u$. If $u^{-1}wu$ is
minimal the result follows so we may assume that this is not the
case. If $u=s\circ t$, where $s$ is non-trivial and commutes with $w$,
then $w^u=w^t$, so the result follows by induction. Therefore we may
assume that $u$ has no such left divisor. From the
Transformation Lemma (see \cite[Lemma 2.5]{BEKR}) it follows  that if
cancellation occurs in the process of reducing $w^u$ to a minimal
form then the first such cancellation must involve a letter from $w$
cancelling with a letter from $u^{\pm 1}$. That is, there exists $y\in
X^{\pm 1}$ such that $u=y\circ v$ and either $w=y\circ w_0$ or $w=w_1\circ
y^{-1}$. In this case $w^u=(w_0\circ y)^v$ or $(y^{-1}\circ
w_1)^v$. In both cases the result follows by induction.
\end{proof}

If $w$ is a minimal form of an element $g\in G$ and $w$ is cyclically minimal
then we say that $g$ is a cyclically minimal element of $G$. It follows from 
Lemma \ref{cycgeod} that if $g$ is a cyclically minimal element of $G$ then all
minimal forms of $g$ are cyclically minimal elements of $(X\cup X^{-1})^*$.
For a cyclically minimal element $g\in G$ we define $[g]=\{h\in G|h=vu, \textrm{
for some } u,v \textrm{ such that } g=u\circ v\}$.

From Proposition 3.9 of \cite{DK}, if $g\in G$ then there exist $u,
w\in (X\cup X^{-1})^*$, with $w$ cyclically minimal, such that
$g=u^{-1}\circ w \circ u$. Thus $g$ is cyclically minimal if and only if
$u=1$. Observe that if $g$ is cyclically minimal then
$l(g^n)=nl(g)$. Therefore partially commutative groups are torsion free.

In Section \ref{sec:cancent} the following lemma and its corollary will
 be useful.
\begin{lem}\label{lem:cmin}
Let $w$ be a cyclically minimal element of $G$ and let $g\in G$. Then
$w^g=u^{-1}\circ v \circ u$ and $g=a\circ b\circ u$, where $v$ is cyclically
minimal, $[v]=[w]$, 
$\a(a)\subseteq\a(w)$, $w^a=v$ and $b\in C(\a(w))$. 
\end{lem}
\begin{proof}
If $l(g)=0$ then the result clearly holds. Assume inductively that $l(g)>1$ and that
the lemma holds for all $g$ of shorter length. If $w^g=g^{-1}\circ w\circ g$ then
the result holds with $u=g$, $w=v$ and $a=b=1$. Therefore we may assume that
$w^g\neq g^{-1}\circ w\circ g$. From \cite[Lemma 2.2]{BEKR} either
{\em (i)} $gcd_l(g,w)\neq 1$; or {\em (ii)} $\gcd_l(g,w^{-1})\neq 1$; or 
{\em (iii)} $D_l(g)\cap C(\a(w))\neq 1$. 

In case  {\em (i)} suppose $\gcd_l(g,w)=c\neq 1$. Then $g=c\circ g_1$,
$w=c\circ w_1$, $w^g=(w_1c)^{g_1}$ and $l(g_1)<l(g)$. As $w$ is cyclically minimal
Lemma \ref{cycgeod} implies that $w_1\circ c$ is also cyclically minimal and by induction
$(w_1c)^{g_1}=u^{-1}\circ v \circ u$, with $[v]=[w_1\circ c]=[w]$, 
$g_1=a^\prime \circ b\circ u$, $\a(a^\prime)\subseteq\a(w_1c)=\a(w)$, 
$(w_1c)^{a^\prime}=v$ and $b\in C(\a(w))$. 
Then $g=c\circ a^\prime\circ b\circ u$, $\a(c\circ a^\prime) \subseteq \a(w)$ and  
$v=w^{c\circ a^\prime}$, so the result holds with $a=c\circ a^\prime$.
A similar argument works in case {\em (ii)}.

Assume then 
that case {\em (iii)} holds.
Let $c\in D_l(g)\cap C(\a(w))$, $c\neq 1$. Then $g=c\circ g_1$, $l(g_1)<l(g)$ and 
$w^g=w^{g_1}$. By induction $w^{g_1}=u^{-1}\circ v\circ u$, where $v$ is
cyclically minimal, $[v]=[w]$, 
$g_1=a\circ b^\prime \circ u$, $\a(a)\subseteq\a(w)$, $w^a=v$ and 
$b^\prime\in C(\a(w))$. As $c\in C(\a(w))$ and $\a(a)\subseteq \a(w)$, we have
 $[c,\a(a)]=1$ and so  
$g=c\circ a\circ b^\prime \circ u= a \circ c\circ b^\prime \circ u$. Therefore
the lemma holds with $b=c\circ b^\prime$.  
\end{proof}
\begin{cor}\label{cor:conj}
Let $w,g$ be (minimal forms of) elements of $G$ and 
$w=u^{-1}\circ v\circ u$, where $v$ is cyclically minimal. Then
there exist minimal forms $a$, $b$, $c$, $d_1$, $d_2$ and 
$e$ such that $g=a\circ b\circ c\circ d_2$, $u=d_1\circ a^{-1}$, $d=d_1\circ d_2$,
$w^g=d^{-1}\circ e\circ d$, $[e]=[v]$, $e=v^b$, $\a(b)\subseteq \a(v)$ and 
$[\a(b\circ c),\a(d_1)]=[\a(c),\a(v)]=1$.
\end{cor}
\begin{proof}
We have $w^g=v^{ug}$. Reducing $ug$ to a minimal form we have 
$u=u_1\circ g_1^{-1}$, $g=g_1\circ g_2$ and $ug=u_1\circ g_2$, for appropriate
$u_1,g_1,g_2$. Then $w^g=v^{u_1\circ g_2}$ and using Lemma \ref{lem:cmin}
$w^g=d^{-1}\circ e\circ d$, $u_1\circ g_2=b\circ c\circ d$,
where $e$ is cyclically minimal, $[e]=[v]$,
$c\in C(\a(v))$, $\a(b)\subseteq \a(v)$ and  $e=v^b$. It follows from Lemma
\ref{lem:fl} that $u_1=b_1\circ c_1\circ d_1$, $g_2=b_2\circ c_2 \circ d_2$, where
$b=b_1\circ b_2$, $c=c_1\circ c_2$, $d=d_1\circ d_2$ and 
$[\a(b_2),\a(c_1\circ d_1)]=[\a(c_2),\a(d_1)]=1$. Then 
$\a(b)\subseteq \a(v)$, and $v^b=e\in [v]$ implies that $\a(b_i)\subseteq \a(v)$ and 
$v^{b_1}\in [v]$. Also $c\in C(\a(v))$ implies $c_i\in C(\a(v))$. 
Now
$u=u_1\circ g_1^{-1}=b_1\circ c_1\circ d_1\circ g_1^{-1}$ and 
$w=u^{-1}\circ v\circ u$, $v^{b_1}\in [v]$ imply $b_1=1$. Then $c_1\in C(\a(v))$ implies
$c_1=1$. Hence $u=d_1\circ g^{-1}_1$, $g=g_1\circ b\circ c\circ d_2$ and 
the result follows with $g_1=a$.
\end{proof}

We shall use one of the normal forms for elements of $G$ defined in
\cite{BEKR}. In order to describe this we introduce the {\em
  non-commutation} graph $\Phi=\Phi(G)$ of $G$. This is the undirected, simple
graph with vertex set $X$ and an edge joining vertices $x$ and $y$ if
and only if $[x,y]\notin R$.
Also, for $Y\subseteq X$, 
denote by $\Phi(Y)$ the full subgraph of $\Phi$ with vertex set
$Y$; so $\Phi(Y)=\Phi(G(Y))$. Extend this notation to arbitrary subsets $S$ of $G$ by setting
$\Phi(S)=\Phi(\a(S))$, where $\a(S)=\cup_{s\in S}\a(s)$. If $S=\{g\}$, a
set with one element, we write $\Phi(g)$ for $\Phi(S)$.
Suppose that the connected components of
$\Phi$ are $\Phi_1$, ..., $\Phi_k$ and let $X_s$ be the vertex set of
$\Phi_s$; so $X=\cup_{s=1}^k X_s$.
 Then $\Phi_s=\Phi(X_s)$ and
$G=G(X_1)\times \cdots \times G(X_k)$.

Let $g$ be a cyclically minimal element of $G$. If $\Phi(g)$ is
connected
then we say that $g$ is a {\em
  block}. Otherwise, suppose that $\Phi(g)$ has connected components
$\Phi_1,\ldots,\Phi_k$. Then $g$ has a minimal form $w=w_1\cdots
w_k$, where $w_s$ is a block with $\Phi(w_s)=\Phi_s$, and consequently
$[w_s,w_t]=1$, for all $s,t\in \{1,\ldots ,k\}$. We call $w_1\cdots
w_k$ a {\em block decomposition} of $g$. Note that it follows from the
Transformation Lemma that if $g$ has another block form $u=u_1\cdots
u_l$ then $k=l$ (as $\Phi(u)=\Phi(g)=\Phi(w)$) and after reordering the
$u_s$'s if necessary $w_s=u_s$, for $s=1,\ldots ,k$. Now let
$g$ be an arbitrary element of $G$ and write $g=u^{-1}\circ w \circ
u$, with $w$ cyclically minimal. Let $w_1\cdots w_n$ be a block
decomposition of $w$. Then $g=w_1^u\cdots w_n^u$ and we say
that $v_i=w_i^u$ is a {\em block} of $g$ and that $g$ has
{\em block decomposition} $v_1\cdots v_n$. In this case a straightforward induction
on $l(u)$ shows that, since $\Phi(w_i)$ is connected, $\Phi(v_i)$ is connected,
for all $i$.

Following \cite{DK}, given $g\in G$  we define $a\in G$ to be a {\em root}
of $g$ if the following conditions hold.
\be[(i)]
\item $a^m=g$, for some $m\ge 1$ and
\item\label{root2} if $b^n=g$, for some $b\in G$ and $n\ge 1$, then $n|m$.
\ee
As shown in \cite{DK} every element of $G$ has a unique root, which we
denote $r(g)$. If $g$ is its own root then we say that $g$ is a {\em
  root element}.
Put another way, $g$ is a root element if $g$ is not a proper power of
any element of $G$. (In fact in \cite{DK} the definition of root is
made in a slightly wider context and also
requires that, in (\ref{root2}), $b=a^{m/n}$. However it is also shown
there that, if $G$ is a partially commutative group and
if $g,h\in G$ with $g^n=h^n$ then $g=h$; so our definition is
equivalent for partially commutative groups.)

If $w\in G$ define $A(w)=G(Y)$, where $Y$ is the set of elements
of $X\backslash \a(w)$ which commute with every element of $\a(w)$. If
$S\subseteq G$ then we define $A(S)=\cap_{w\in S}A(w)$.
Let $w$ be a cyclically minimal element of $G$
with block decomposition $w=w_1\cdots w_k$ and
let $v_i=r(w_i)$. Then,
from \cite[Theorem 3.10]{DK},
\begin{equation}\label{eq:centraliser}
C(w)=\langle v_1\rangle \times \cdots \times \langle v_k \rangle\times A(w).
\end{equation}
It follows  that
\begin{equation}\label{eq:cent_int}
C(w)=\displaystyle{\cap_{i=1}^k} C(v_i).
\end{equation}
Furthermore, if $u=z^{-1}\circ w \circ z$ then $C(u)=[C(w)]^z$ so
\begin{equation}\label{eq:centraliser_conj}
C(u)=
\langle v_1^z\rangle \times \cdots \times \langle v_k^z \rangle\times
[A(w)]^z,
\end{equation}
but note that in general $[A(w)]^z\neq A(u)$.

We shall also need the following lemma, which 
is based on Proposition 5.7 of \cite{BEKR}, and its corollary.
\begin{lem}\label{conj_block}
Let $u$ and $v$ be cyclically minimal with block decompositions
$u=u_1\cdots u_k$ and $v=v_1\cdots v_l$. If $u=v^g$, for some $g\in G$
then $k=l$ and, after reordering of the $v_s$'s if necessary,
$u_s=v_s^g$ and some minimal form of $u_s$ is equal to  a cyclic permutation of $v_s$, for
$s=1,\ldots ,k$.
\end{lem}
\begin{proof}
It follows from the proof of Proposition 5.7 of \cite{BEKR} that $k=l$
and that, after an appropriate reindexation, $u_s=v_s^a$,  for some
$a\in G$ and for all $s$. Then $u=v^a=v^g$ which implies that $v=v^{ag^{-1}}$, so
$ag^{-1}\in C(v)$. This implies that $g=za$, for some $z\in C(v)$, and
so, from \eqref{eq:cent_int}, $z\in C(v_s)$, for $s=1,\ldots
,k$. Thus $u_s=v_s^a=v_s^g$, for all $s$, as claimed. The final
assertion follows since $u_s$ is minimal.
\end{proof}

\begin{cor}\label{conjcent}
Let $u$ and $v$ be cyclically minimal, root, block elements and let $g$
be such that $u^g=g^{-1}\circ u \circ g$. If
$u$ and $u^g$ are in $C(v)$  then $g\in C(v)$.
\end{cor}
\begin{proof}
From \eqref{eq:centraliser} and the fact that $u$ and $v$ are block,
root elements
and $u\in C(v)=\la v\ra\times A(v)$ we have $u\in \la v\ra$ or $u\in
A(v)$. Thus either $u=v^{\pm 1}$, and without loss of generality
$u=v$, or
$u\in A(v)$.

Suppose first that $u\notin A(v)$, so $u=v$. As $u^g$ is a block, $\Phi(u^g)$ is
connected so
$u^g\in C(v)$ implies $u^g\in \la v\ra$
or $u^g\in A(v)$.
Since $u^g$ is minimal as written, 
if $u^g\in A(v)$ then 
$u\in A(v)$, a contradiction.
As $v$ is cyclically minimal, if $u^g\in \la v\ra$,
and the given hypotheses hold, then $g=1$; so the result follows.

Now suppose that $u\in A(v)$. Then, since it's a block,  $g^{-1}\circ u
\circ g\in A(v)$, so $g\in A(v)$ 
and the
result follows.
 \end{proof}
\section{Centraliser dimension and canonical centralisers}\label{sec:cancent}
In this section we prove that the height of the centraliser
lattice $\fC(G)$ of $G$ coincides with the height of a sublattice: namely
the lattice of canonical
centralisers $\fC(X;G)$. First we must show that $\fC(X;G)$ is indeed a sublattice
of $\fC(G)$ and to do so we need to describe the structure of canonical
centralisers.

Let $Y$ be a subset of $X$ and define the set of generators {\em
  orthogonal} to $Y$ to be
$Y^\bot=\{x\in X|[x,y]=1,
\textrm{ for all } y\in Y\}$. Then $G(Y^\bot)\subseteq C(Y)$.
On the other hand,
if $z\notin G(Y^\bot)$ then there exists $x\in \a(z)$ such that
$x\notin Y^\bot$. Therefore there is $y\in Y$ such that $[x,y]\neq 1$,
that is $x\notin C(y)$.
Now
from \eqref{eq:centraliser}, $C(y)=\la y\ra \times A(y)$, so $x\notin
C(y)$ implies $z\notin C(y)$. It follows that $C(Y)=G(Y^\bot)$, a
parabolic subgroup. To rephrase this in terms of the notation of
Section \ref{sec:prelim} let $Y_0=Y\cap Y^\bot =\{y_1,\ldots ,y_n\}$ and
$Z=Y^\bot\backslash Y$, so $Y^\bot=Y_0\cup Z$. Then $y\in Y_0$ implies
$[y,u]=1$, for all $u\in Y^\bot$,
and  moreover $G(Z)=A(Y)$. Therefore we have;
\begin{equation}\label{eq:cancent}
\textrm{ if } Y\subseteq X \textrm{ then } C(Y)=G(Y^\bot)=\la y_1\ra
\times \cdots \times \la y_n\ra \times A(Y),
\end{equation}
where $\{y_1,\ldots ,y_n\}= Y^\bot\cap Y$.
\begin{lem}\label{lem:cancent}
The set of canonical centralisers $\fC(X;G)$ is a sublattice
of $\fC(G)$.
\end{lem}
\begin{proof}
It is clear that the intersection of two canonical centralisers
is a canonical centraliser. Now let $C_1=C(Y_1)$ and $C_2=C(Y_2)$
be canonical centralisers,
for some subsets $Y_1$ and $Y_2$ of $X$. Then, in the lattice $\fC(G)$ we
have $C_1\vee C_2=C$, where $C$ is the intersection of all centralisers containing
$C_1$ and $C_2$. Then $C=C(W)$, for some subset $W\subseteq G$, and we set
$Z=\cup_{w\in W} \a(w)$, so $Z\subseteq X$. Clearly $C(Z)\subseteq
C(W)$.
From \eqref{eq:cancent}, $C(Y_1)=G(Y_1^\bot)$. If $a\in Y^\bot_1$ then
$a\in C_1\subseteq C(W)$, so $a\in C(w)$, for all $w\in W$. As $a\in
X$ it follows that $a\in C(\a(w))$, for all $w\in W$. Therefore
$Y^\bot_1\subseteq C(Z)$ and so $C_1\subseteq C(Z)$.
Similarly
$C_2\subseteq C(Z)$ so we have $C=C_1\vee C_2\subseteq C(Z)\subseteq C$. Hence
$C=C(Z)$, a canonical centraliser. That is $C_1\vee C_2\in \fC(X;G)$
and this completes the proof of the lemma.
\end{proof}

Expression \eqref{eq:centraliser} for the centraliser of cyclically
minimal elements gives the structure of  the centraliser of certain
sets of more than one element as follows.
\begin{lem}\label{k-std-form}
Let $w_1,\ldots, w_{k-1},u$ be cyclically minimal root, block elements of
$G$. Let $w_k=g^{-1}\circ u\circ g$, for some $g\in G$ such that
either $g=1$ or $g\notin C(w_1,\ldots ,w_{k-1})$. Suppose
that
\begin{equation}\label{eq:sat}
C(w_1)\sdc C(w_1,w_2)\sdc \cdots \sdc C(w_1,\ldots ,w_k)
\end{equation}
is a saturated chain of centralisers.
Then, there exist integers $m\ge 0$ and
$i_1,\ldots ,i_m$, with
$1\le i_1< \cdots <i_m\le k-1$, such that
\be[(i)]
\item
if $g=1$ and $w_k\in
C(w_1,\ldots ,w_{k-1})$ then
\begin{equation}\label{eq:centraliser_many1}
C(w_1,\ldots ,w_k)=\la w_{i_1}\ra \times \cdots \times \la w_{i_m}\ra
\times \la w_k\ra
\times A(w_1,\ldots ,w_k),
\end{equation}
\item and otherwise
\begin{equation}\label{eq:centraliser_many2}
C(w_1,\ldots ,w_k)=\la w_{i_1}\ra \times \cdots \times \la w_{i_m}\ra
\times (A(w_1,\ldots ,w_{k-1})\cap A^g(u)).
\end{equation}
\ee
\end{lem}
\begin{proof}
When $k=1$ this is \eqref{eq:centraliser} or
\eqref{eq:centraliser_conj}. Write $C_i=C(w_1,\ldots, w_i)$ and
$A_i=A(w_1,\ldots ,w_i)$. Assume inductively that
\[
C_{k-1}=\la w_{i_1}\ra \times \cdots \times \la w_{i_n}\ra
\times A_{k-1},
\]
where $n\ge 0$ and $1\le i_1< \cdots <i_n\le k-1$.
(Here $w_{i_n}=w_{k-1}$ if and only if $w_{k-1}\in C_{k-2}$.)
We have
\begin{equation}\label{eq:indhyp}
C_k=C_{k-1}\cap (\la w_k\ra\times
A^g(u)).
\end{equation}
If $g=1$ and
$w_k=u\in A_{k-1}$ then $w_{i_1},\ldots,w_{i_n}\in C(w_k)$.
As $w_{i_j}$ is a block this implies that $w_{i_j}\in \la w_k\ra$ or
$w_{i_j}\in A(w_k)$. If $w_{i_j}\in \la w_k\ra$ then
$C(w_{i_j})=C(w_k)$, so \eqref{eq:sat} is not strictly descending,
contrary to hypothesis. Hence $w_{i_j}\in A(w_k)$, for $j=1,\ldots,n$.
A straightforward check, using
\ref{eq:indhyp}, now shows that
\[
C_k=\la w_{i_1}\ra \times \cdots \times \la w_{i_n}\ra
\times \la w_k\ra
\times A_k,
\]
as required.

If $g=1$ and $w_k\notin A_{k-1}$ then $w_k\in C_k$ implies $w_k\in \la
w_i \ra$, for some $i$ with $1\le i\le k-1$. This implies
$C(w_k)=C(w_i)$, a contradiction. Hence we may assume $g\neq 1$
or $w_k\notin C_{k-1}$.

If $g\neq 1$ and $w_k\in A_{k-1}$ then, since $u^g$ is minimal, we
have $g\in A_{k-1}\subseteq C_{k-1}$, a contradiction.
Suppose then that $w_k\notin A_{k-1}$.
Let $z\in C_k$. Then, from \ref{eq:indhyp},
\[
z=w_{i_1}^{\a_1}\cdots w_{i_n}^{\a_n}\cdot a=w_k^\b\cdot b^g,
\]
where $a\in A_{k-1}$ and $b\in A(u)$.
Write $a=a^{-1}_0\circ a_1\circ a_0$ and $b=b_0^{-1}\circ b_1\circ b_0$, where $a_1$ and $b_1$
are cyclically minimal. Note that we have $a_0,a_1\in A_{k-1}$ and
$b_0,b_1\in A(u)$. Then
\begin{align*}
w_{i_1}^{\a_1}\cdots w_{i_n}^{\a_n}\cdot a_1 & = a_0 (w_k^\b
b_1^{b_0g})a_0^{-1}\\
 & = a_0(u^\b b_1)^{b_0g}a_0^{-1}, \textrm{ since } b_0\in A(u),\\
 & = (u^\b b_1)^{b_0 g a_0^{-1}}.
\end{align*}
Now $z_1=w_{i_1}^{\a_1}\cdots w_{i_n}^{\a_n}\cdot a_1$ and $u^\b b_1$
are cyclically minimal and, if $\b\neq 0$ then  $u^\b$ is a block of
$u^\b b_1$. Therefore,
from Lemma \ref{conj_block}, if $\b\neq 0$ then
$(u^\b)^{b_0ga_0^{-1}}$ is a block of $z_1$. Blocks of $z_1$ are
either of the form $w_{i_j}^{\a_j}$ or are blocks of $a_1$. However
\begin{align}
(u^\b)^{b_0ga_0^{-1}} & = w_{i_j}^{\a_j} \textrm{ implies }\notag\\
(u^\b)^{ga_0^{-1}} & = w_{i_j}^{\a_j} \textrm{ implies }\notag\\
(u^\b)^{g} & =  w_{i_j}^{\a_j}, \textrm{ since } a_0 \in A_{k-1}.\label{al:block}
\end{align}
If $g\neq 1$ then 
$(u^\b)^{g}$ is not cyclically minimal, and so \eqref{al:block} is 
a contradiction. On the other hand if $g=1$ then \eqref{al:block} implies that
$u^\b=w_{i_j}$ so $w_k=u\in C_k$, again a contradiction. Hence,
if $\b\neq 0$ then $(u^\b)^{b_0ga_0^{-1}}$ is a block of $a_1$ and so
belongs to $A_{k-1}$. This implies that $(u^\b)^{g}\in A_{k-1}$ and,
since $u$ is cyclically minimal, that $w_k=u^g\in
A_{k-1}$, another contradiction. Thus $\b=0$ and
$z_1=b_1^{b_0ga_0^{-1}}$. From Lemma \ref{conj_block} again, for
$j=1,\ldots ,n$, either $\a_j=0$ or
\[
w_{i_j}^{\a_j}=b_{s_j}^{b_0ga_0^{-1}} \implies w_{i_j}^{\a_j}=b_{s_j}^{b_0g},
\]
for some block $b_{s_j}$ of $b_1$. Thus, if $\a_j\neq 0$ then  $w_{i_j}\in A(u)^g\subseteq
C(w_k)$ and, since by definition $w_{i_j}\in C_{k-1}$, we have
$w_{i_j}\in C_k$ or $\a_j=0$, for $j=1,\ldots,n$. Also
$w_{i_1}^{\a_1}\cdots w_{i_n}^{\a_n}\in A^g(u)$ and so  $a\in
A_{k-1}\cap A^g(u)$. This means that $z\in T$, where
\[
T=\la w_{j_1}\ra \times \cdots \times \la w_{j_m}\ra
\times (A_{k-1}\cap A^g(u)),
\]
with $\{j_1,\ldots,j_m\}$ the subset of $\{1,\ldots ,k-1\}$ such
that
$w_{j_l}\in C_k$. Hence $C_k\subseteq T$.
The opposite
inclusion clearly holds and  the result follows.
\end{proof}

The centraliser dimension of a partially commutative group $G$ is by definition
the height of its centraliser lattice $\fC(G)$. The main theorem of this paper,
which follows, shows that this height is the same as the height of the sublattice
consisting of centralisers of canonical generators. This means that the centraliser
dimension may be calculated directly from the commutation graph of $G$.
\begin{thm}\label{thm:pcdimgen}
Let $G$ be a partially commutative group. Then centraliser dimension
of $G$  is equal to the height of the lattice $\fC(X;G)$ of canonical
centralisers of $G$.
\end{thm}
\begin{proof}
As pointed out above, $G$ has finite centraliser dimension as it is a linear group.
Let $d$ be the centraliser dimension of $G$. We must show that
there exist elements $x_1, \ldots ,x_d\in X$ such that
\[
G>C(x_1)\sdc\cdots \sdc C(x_1, \ldots ,x_d)=Z(G)
\]
is  a centraliser chain in $G$. Let
\begin{equation}\label{eq:pcmax}
  G=C(1)\sdc C(w_1) \sdc \cdots \sdc C(w_1,\ldots ,w_d)=Z(G)
\end{equation}
be a maximal centraliser chain in $G$. For $i=1,\ldots, d$, set
$C_i=C(w_1,\ldots ,w_i)$
and $A_i=A(w_1,\ldots ,w_i)$.

\noindent{\bf Step 1.} 
For $i=1,\ldots ,d$
we may write $w_i=g_i^{-1}\circ u_i\circ g_i$, with $u_i$ cyclically
minimal. We shall show in this step that we may assume that $u_i$ is
also a block, root element.

Suppose
that the block decomposition of $u_j$ has exactly one block, for $j=1,\ldots ,i-1$, but
the block decomposition of $u_i$ is $u_i=v_1\cdot \cdots \cdot v_k$, with $k>1$.
As $C_{i-1}\neac C(w_i)$ it follows that $C_{i-1}\neac g_i^{-1}C(v_t)g_i$, for some
$t$. As \eqref{eq:pcmax} is maximal
and $g_i^{-1}C(v_t)g_i\sdc C(w_i)$
it follows that $g_i^{-1}C(v_t)g_i\cap C_{i-1}=C(w_i)\cap C_{i-1}$.
Thus we may assume that $u_i$ is a block, for $i=1,\ldots , d$. We may also assume that
$u_i$ is a root element, since $C(u)=C(v)$, if $v$ is the root of a block $u$.
Thus
\[
C(w_i)=\langle u_i \rangle ^{g_i}\times A^{g_i}(u_i).
\]
\noindent{\bf Step 2.} 
Now we shall show that we may in addition assume that $w_i$ is cyclically minimal, for all $i$.
Assume that $g_j=1$, for $j=1,\ldots ,i-1$, but that $g_i\neq 1$, that is, $w_i$ is not cyclically
minimal.  First consider the case where 
$D_r(g_i)\cap  C_{i-1}\neq 1$, so $g_i=a\circ b$, where $b\in C_{i-1}$.
Then, writing $h=b^{-1}$, we have
$C_j^h=C_j$, for $j=1,\ldots ,i-1$. Hence
\[
G\sdc C_1\sdc \cdots \sdc C_{i-1}=C_{i-1}^h\sdc C_i^h \sdc C_d^h=Z(G)
\]
is a maximal centraliser chain. Since $C_i^h=C(w_1,\ldots ,w_{i-1},u_i^{a})$ 
we may replace
$w_i$ with $u_i^a$, in this case. Therefore we may assume that $D_r(g_i)\cap  C_{i-1}=1$.
In particular $g_i\notin C_{i-1}$.
Since $w_1, \ldots, w_{i-1}$ are cyclically minimal, from Lemma
\ref{k-std-form} we have
\[
C_j=\langle w_{j,1}\rangle \times \cdots \times \langle w_{j,k_j} \rangle \times A_j,
\]
where the $w_{j,i}$ are in $\{w_1,\ldots ,w_j\}$,
for all $j$ with $1\le j\le i-1$. In particular we can write
\[
C_{i-1}=\langle w_{i_1}\rangle \times \cdots \times \langle w_{i_m} \rangle \times A,
\]
where $1\le i_1<\cdots <i_m\le i-1$ and $A=A_{i-1}$. If $w_i\in A$ then, since
$u_i^{g_i}$ is minimal we have $g_i\in A\subseteq C_{i-1}$, a
contradiction. Hence $w_i\notin A$ and from Lemma \ref{k-std-form}, after
reordering $i_1,\ldots i_m$, if necessary, we have
\begin{equation}\label{eq:C_i}
C_i=\langle w_{i_1}\rangle \times \cdots \times \langle w_{i_k} \rangle \times (A \cap A^{g_i}(u_i)),
\end{equation}
where $w_{i_1},\ldots , w_{i_k}\in A^{g_i}(u_i)$, and
$w_{i_{k+1}},\ldots ,w_{i_m}\notin A^{g_i}(u_i)$,
for some $k\le m$.

Next we show that $A\cap A^{g_i}(u_i)\subseteq A(u_i)\cap C(g_i)$. Suppose that
$h\in  A\cap A^{g_i}(u_i)$, so there is $f\in A(u_i)$ such that $h=f^{g_i}$.
Write $f=k^{-1}\circ l \circ k$, where $l$ is cyclically minimal. From
Corollary \ref{cor:conj}, there are elements $a$, $b$, $c$, $d_1$, $d_2$, $e$,
such that 
$g_i=a\circ b\circ c\circ d_2$, $k=d_1\circ a^{-1}$, $d=d_1\circ d_2$, 
$h=f^{g_i}=d^{-1}\circ e\circ d$, $[e]=[l]$, $e=l^b$, $\a(b)\subseteq \a(e)$ and 
$1=[\a(b\circ c),\a(d_1)]=[\a(c),\a(l)]$. As $h\in A$ we have
$d_1,d_2\in A$ and 
$D_r(g_i)\cap C_{i-1}=1$, 
so $d_2=1$ and $d=d_1$. Also $f\in A(u_i)$
so $a,d,l \in A(u_i)$ and $\a(b)\subseteq \a(e)$, $[e]=[l]$, implies
$e,b\in A(u_i)$. Since $D_l(g_i)\cap C(u_{i})=1$, this implies that
$a=b=1$ so $g_i=c$, $d=k$, $e=l$ and $1=[\a(g_i),\a(d)]=[\a(g_i),\a(e)]$. 
Hence $h=f\in A(u_i)$ and $h\in C(g_i)$, so
 $A\cap A^{g_i}(u_i)\subseteq A(u_i)\cap C(g_i)$, as claimed.

Hence $A\cap A^{g_i}(u_i)\subseteq A\cap A(u_i)\cap
C(g_i)$. Conversely if $c\in A\cap A(u_i)\cap C(g_i)$ then $c\in
A\cap A^{g_i}(u_i)$. Therefore
\begin{equation}\label{eq:Aform}
A\cap A^{g_i}(u_i)=A\cap A(u_i)\cap C(g_i).
\end{equation}
From \eqref{eq:C_i} and \eqref{eq:Aform} we have
\begin{equation}\label{eq:C_i_x}
C_i=\langle w_{i_1}\rangle \times \cdots \times \langle w_{i_k}
\rangle \times (A \cap A(u_i)\cap C(g_i)).
\end{equation}

Note that  
$\langle w_{i_1}\rangle \times \cdots \times \langle w_{i_k} \rangle 
\subseteq  A\cap A^{g_i}(u_i)$
implies that
$\langle w_{i_1}\rangle \times \cdots \times \langle w_{i_k} \rangle 
\subseteq A\cap A(u_i)$.
Therefore
\begin{equation*}
C_{i-1}\cap C(u_i) \edc \langle w_{i_1}\rangle \times \cdots \times \langle w_{i_k}\rangle
\times (\langle u_i \rangle \cap  A)\times ( A\cap  A(u_i))\edc C_i,
\end{equation*}
from \eqref{eq:C_i_x}.
If $C_{i-1}\subseteq C(u_i)\cap C(g_i)$ then $C_{i-1}\subseteq
C(u_i^{g_i})=C(w_i)$, which cannot happen as
the chain \eqref{eq:pcmax} is maximal.
If $C_{i-1}\nsubseteq C(u_i)$ then $C_{i-1}\cap C(u_i)\sac C_{i-1}$
and we may replace $w_i$ by $u_i$, since by maximality of
\eqref{eq:pcmax}, $C_i=C_{i-1}\cap C(u_i)$ in
this case. Therefore we may assume that $C_{i-1}\subseteq C(u_i)$ but
that $C_{i-1}\nsubseteq C(g_i)$.
As $C_{i-1}\subseteq C(u_i)$ we have $w_{i_j}\in C(u_i)$, for
$j=1,\ldots k$. Hence $w_i=u_i^{g_i}\in C(w_{i_j})$ and $u_i\in
C(w_{i_j})$,
for $j=1,\ldots ,k$. From Corollary \ref{conjcent}, $g_i\in
C(w_{i_j})$ and so, from \eqref{eq:C_i_x}, $C_i\subseteq C(g_i)$.
Since \eqref{eq:pcmax} is maximal and we have $C_{i-1}\sdc C_{i-1}\cap
C(g_i)\edc C_i$ it follows that $C_{i-1}\cap C(g_i)=C_i$. Hence we may
(temporarily) replace $w_i$ by $g_i$ and, as in Step 1,
write $g_i=h_0^{-1}\circ h_1\circ h_0$, with $h_1$
cyclically minimal and then replace $g_i$ by $h_0^{-1}h_2h_0$, where
$h_2$ is the root of a block of $h_1$. Thus we have replaced $w_i$ with
$h_0^{-1}h_2h_0$ and $l(h_0^{-1}h_2h_0)<l(w_i)$. We now repeat this
process until, eventually, we replace $w_i$ by a cyclically minimal,
root, block element. By induction we may therefore assume that $w_i$
is cyclically minimal, for $i=1,\ldots ,d$.\\

\noindent{\bf Step 3.}
Finally we show that we may replace each $w_i$ with an element of $X$.
As $C_d=Z(G)$ we have, using Lemma \ref{k-std-form},
$C_d\eac  A(w_1,\ldots, w_d)$, otherwise $C_i=C_{i+1}$, for some $i$.
Hence, for each $i=1,\ldots , d$ there exists a minimal integer $j=d(i)$ such that $w_i\notin C_j$.
Suppose that for some $i\ge 1$ we have $j=d(i)\ge i+1$ and $l(w_i)>1$. Let $r_i$ be a
cyclically minimal block
such that $\az(r_i)=\az(w_i)$, $r_i$ is a root element
and $r_i$ does not commute with $w_i$: for example a cyclic permutation
of $w_i$. Then
\[
C(w_1,\ldots,w_{j-1},r_i)=C_{j-1}\cap  A(w_i)\edc C_j.
\]
Since \eqref{eq:pcmax} is maximal and $w_i\in C_{j-1}\backslash  A(w_i)$ we have
$C_{j-1}\cap  A(w_i)=C_j$. Therefore we may replace $w_j$ with
$r_i$. If $j=i+1$ this is all we need do.

Otherwise, if $j>i+1$, define $C^{(i+1)}_i=G$ and,  for $t=i+1,\ldots ,j-1$,
define $C_t^{(i+1)}=C(w_{i+1},\ldots ,w_t)$ and
$C^\prime_t=C_{i-1}\cap C_t^{(i+1)}\cap  A(w_i)$.
Then
$C^\prime_t=C_{i-1}\cap C(w_i,r_i)\cap C_t^{(i+1)}$ is a centraliser of $t+1$ elements of $G$,
for $t=i,\ldots , j-1$. Also
$C^\prime_{j-1}=C_j$, by definition of $r_i$,  and
$C_t=\langle w_i \rangle \times C^\prime_t$, for $t=i,\ldots , j-1$.
As \eqref{eq:pcmax} is a saturated centraliser chain so is
 \[
C^\prime_i\sdc \cdots \sdc C^\prime_{j-1}.
\]
Hence
\[
G\sdc
C_1\sdc \cdots \sdc C_i \sdc C^\prime_i\sdc \cdots \sdc C^\prime_{j-1}
=C_j\sdc \cdots \sdc C_d=Z(G)
\]
is a maximal centraliser chain. We may therefore replace the sequence
$w_1,\ldots ,w_d$ with the sequence
$w_1,\ldots, w_i,r_i,w_{i+1},\ldots ,w_{j-1},w_{j+1},\ldots ,w_d$ to give a new chain
\eqref{eq:pcmax}. In the new chain $d(i)=i+1$. Note that $w_i\in  A_{i-1}$ so that
if $h<i$ then $d(h)\neq i$. Hence, beginning with $i=1$, we may repeat this process for successive
values of $i$ until,
for all $i$, either $l(w_i)=1$ or $d(i)\le i+1$.

Now suppose that $l(w_i)=1$ and that $j=d(i)\ge i+1$.
Write $w_i=x$ for notational convenience. Let $y$ be a letter of $\az(w_j)$ which does not commute
with $x$. (Such a $y$ must exist as $d(i)=j$.) Then $ A(y)\edc  A(w_j)$ so
\[
C_{j-1}\sdc C_{j-1}\cap C(y) \edc C_{j-1}\cap C(w_j)=C_j.
\]
By maximality of \eqref{eq:pcmax} we may replace $w_j$ with $y$. Repeating this for all $i$ such
that $l(w_i)=1$ and that $j=d(i)\ge i+1$, we obtain a new sequence $w_1, \ldots ,w_d$ with the property
that either
\be[(i)]
\item\label{step3i} $l(w_i)=1$, or
\item\label{step3ii} $d(i)=d(i+1)=i+1$, in which case we call $w_i$ and $w_{i+1}$ a
{\em  matched pair}, or
\item $d(i)\le i$, $l(w_i)>1$ and $w_i$ is not one of a matched pair,
\ee
for all $i=1,\ldots,d$.

Suppose that $i$ is such that neither \eqref{step3i} or
\eqref{step3ii} holds. Then $w_i\notin C_{i-1}$,
since this would imply $w_i\in A_{i-1}$ and so, from Lemma 
\ref{k-std-form}, $d(i)\ge i+1$.
Therefore, using Lemma  
\ref{k-std-form} again, $C_i=C_{i-1}\cap A(w_i)$. As $C_i<C_{i-1}$ and
$A(w_i)=A(\a(w_i))$ there exists $y\in\a(w_i)$ such that
$A(y)\nsupseteq C_{i-1}$. Since $A(w_i)\le A(y)$ we have
\[
C_{i-1}\sdc C_{i-1}\cap A(y)\edc C_{i-1}\cap A(w_i)=C_i.
\]
Hence we may replace $w_i$ with $y$. This results in one fewer
occurrence of an index $i$ which does not satisfy \eqref{step3i} or
\eqref{step3ii}. Repeating this process we therefore obtain a sequence
$w_1,\ldots, w_d$ such that \eqref{step3i} or
\eqref{step3ii} holds, for all $i$.

Consider now $i$ such that $l(w_i)>1$ and $d(i)=i+1$. Then as above we may assume $w_{i+1}=r_i$. Choose letters
$x$, $y\in w_i$ which do not commute. These exist as $w_i$ is a block of length more than one. As
$\az(w_i)\subseteq  A_{i-1}$ we have $x,y\in  A_{i-1}$, so $x,y\in C_{i-1}$.
Define $C^\prime_i=C_{i-1}\cap C(x)$ and $C^\prime_{i+1}=C_{i-1}\cap C(x,y)$.
Then $C(x,y)= A(x,y)$ so $C^\prime_{i+1}=C_{i-1}\cap  A(x,y)\edc C_{i-1}\cap  A(w_i)=C_{i+1}$.
As $y\notin C(x)$ and $x\notin C(x,y)$ we have
\[
\cdots\sdc C_{i-1}\sdc C^\prime_i\sdc C^\prime_{i+1}\edc C_{i+1}\sdc \cdots
\]
and, by maximality of \eqref{eq:pcmax}, $C^\prime_{i+1}=C_{i+1}.$ We may now replace $w_i,w_{i+1}$ by
$x,y$. Making such replacements over all such pairs $w_i,w_{i+1}$ gives the required result.
\end{proof}

\section{Extensions of Partially Commutative Groups}\label{sec:extpc}

In this section we investigate the effect on centraliser dimension of
adding or removing a vertex (and its incident edges) from the
commutation graph of a partially commutative group. We shall see that
addition
of a vertex can either leave the centraliser dimension unaltered or
increase it by one or two.

We call a sequence of generators $x_1,\ldots ,x_l$ of a partially commutative group $G$ a
{\em parameter system} for $G$ if
$$G=C_0\sdc C_1 \sdc \cdots \sdc C_n\supseteq Z(G),$$
where $C_i=C(x_1,\ldots ,x_i)$ is a centraliser
chain in $G$. We say that the parameter system is {\em maximal} or
{\em saturated} if the centraliser chain is maximal or saturated,
respectively. If $C_n=Z(G)$ then we say that the parameter system {\em ends} in
$Z(G)$.

Let $G$ be a partially commutative group with commutation graph $\G$
and presentation $\la
X|R\ra$, as before, and fix $x\in X$. Let $X_x=X\backslash \{x\}$, let
$\G_x$ be the graph obtained from $\G$ by deleting $x$ and
$R_x=\{r\in R: r=[a,b], a\neq x,b\neq x\}$. Let $G_x$ be the
partially commutative group with commutation graph $\G_x$; so $G_x$ has presentation $\la X_x|R_x\ra$.
Let
$
Y =Y_x= \{y\in X_x: [y,x]=1 \textrm{ in } G\},
$
let $W=W_x=X_x\backslash Y$
and, as usual,
let $A=A_G(x)=G(Y)$; so $A$ is a canonical parabolic subgroup of both
$G$ and $G_x$.
Since the above presentation for $G$ can be written as
$$
\la X_{x}, x \mid R_x,
\left[x, y \right]= 1 \hbox{ for } y \in
Y \ra
$$
it
is easy to see that $G$
is the $HNN$-extension
$\la G_{x},x\mid x^{-1}ax=a, \forall a\in A \ra$.

Now suppose that $P=x_1,\ldots, x_n$ is a parameter system for
$G_x$ and let $C_i=C_{G_x}(x_1,\ldots,x_i)$.
In order to understand the effect on centraliser dimension of forming
the $HNN$-extension $G$ from $G_x$ we need the following definitions.
If there is some $l$ such that $\{x_1,\ldots ,x_l\}\subseteq Y$
and some $w\in W$ such that  $C_{G_x}(w)\supseteq C_l$ then we say
that
$P$ is {\em locked} at $l$ with key $w$ (with respect to $A$).
As we shall see below, if $P$ is locked at $l$ with key $w$ then the sequence
$x_1,\ldots ,x_l, w,x_{l+1},\ldots ,x_n$ is a parameter system for $G$.

Now let $k$ be maximal such that $C_k\nsubseteq A$.
We say that $P$ is {\em tied} at $k$, of type T\ref{it:t1a},
T\ref{it:t1b} or T\ref{it:tied2}, at $k$ if either
\be[{T}1.]
\item\label{it:tied1}
\be[a.]
\item\label{it:t1a}
$k<n$ and there exists $y\in Y$ such that $y\in
  C_k\backslash C_{k+1}$; or
\item\label{it:t1b}  $k=n$ and there exists $w\in W\cap C_n$; or
\ee
\item\label{it:tied2} $k<n$, $\{x_1,\ldots ,x_k\}\subseteq Y$ and
  $x_{k+1}\in W$.
\ee
Again, if $P$ is tied at $k$ then it turns out that
$x_1,\ldots ,x_k, x,x_{k+1},\ldots ,x_n$ is a parameter system for
$G$.

Given a maximal parameter system $P$ for $G_x$ we shall see that $P$
is a maximal parameter system for $G$ if and only if $P$
is neither locked nor tied.
Moreover we shall
describe below how any maximal parameter
system $Q$ for $G$ induces a parameter system for $G$ and see that if this
parameter system is shorter than $Q$ then it is either locked or tied.
As a result it will follow (see Theorem \ref{cd-ext}) that the centraliser dimension of $G$ is at
most $2$ more than that of $G_x$.

Before using these definitions to state our theorem we illustrate them
with some examples.
\begin{expl}\label{ex:lt}
In each of the following cases $G$
is an $HNN$-extension of $G_x$ with stable letter $x$, as above.
\be
\item\label{it:lt1}
In the group with commutation graph $\G_1$ of Figure \ref{fig:lt1}
we have $Y=\{a,b,c\}$  and $W=\{d,e,f\}$. By inspection $\cd(G)=5$ and
$\cd(G_x)=4$.
Among the maximal parameter systems of $G_x$ are the following.
\be[(i)]
\item $d,f,a,c$  which is tied at $2$ $(x_2=f)$, of type
T\ref{it:t1a}, but not locked. The corresponding centraliser chain
is
$$
G_x\sdc \la d,a,b,c\ra\sdc \la a,c\ra\sdc \la a\ra \sdc 1 .
$$
\item $c,a,b,f$ with centraliser chain
$$
G_x\sdc \la f,d,e\ra\sdc \la f,d\ra\sdc \la d\ra \sdc 1 .
$$
This is tied of type
T\ref{it:tied2} at $3$ $ (x_3=b)$ and is also locked at $3$ with key $d$.
\item  $c,a,f,b$ which is neither locked nor tied.
 The corresponding centraliser chain
is
$$
G_x\sdc \la f,d,e\ra\sdc \la f,d\ra\sdc \la f\ra \sdc 1 .
$$
\ee
\item\label{it:lt2}
In the group with commutation graph $\G_2$ of Figure \ref{fig:lt1}
we have $Y=\{a,b,c\}$  and $W=\{d,e\}$. In this case $\cd(G)=6$ and $\cd(G_x)=4$.
$G_x$ has a maximal parameter system
$b,a,c,e$
with centraliser chain
$$
G_x\sdc \la a, b,c, d\ra\sdc \la a,b,d\ra\sdc \la b,d\ra \sdc \la d\ra
.
$$
This is locked at $3 $ $ (x_3=c)$, with key $d$,  and tied at $4 $ $ (x_4=e)$,
of type T\ref{it:t1b}.
\item\label{it:lt3}
In the group with commutation graph $\G_3$ of Figure \ref{fig:lt1}
we have $Y=\{c,d\}$  and $W=\{a,b\}$. Here $\cd(G)=\cd(G_x)=4$.
$G_x$ has a  maximal parameter system
$c,b,d,a$
with centraliser chain
$$
G_x\sdc \la a, b,c\ra\sdc \la b,c\ra\sdc \la b\ra \sdc 1
.
$$
This is neither locked nor tied.
\item\label{it:lt4} In the group with commutation graph $\G_4$ of Figure \ref{fig:lt2}
we have $Y=\{a,b,c\}$  and $W=\{d,e\}$. Here $\cd(G)=6$ and $\cd(G_x)=4$.
$G_x$ has a  maximal parameter system
$b,a,c,e$
with centraliser chain
$$
G_x\sdc \la a, b,c,d\ra\sdc \la a,b,d\ra\sdc \la b,d\ra \sdc 1
.
$$
This is locked at $3$ $(x_3=c)$, with key $d$, and  tied at $3$, of type
T\ref{it:t1a}.
\item In the group with commutation graph $\G_5$ of Figure \ref{fig:lt2}
we have $Y=\{a,b,c\}$  and $W=\{d,e,f,g,h,i,j,k\}$. Here $\cd(G_x)=5$,
since its commutation graph has degree $3$ and $l,i,k,j,h$ is a parameter system for
$G_x$.
From  Example \ref{ex:lt}.\ref{it:lt4} above $G$ must have a parameter system
of
length $6$ and as $\G_5$ has degree $4$ it follows that $\cd(G)=6$.
We remark that from Example \ref{ex:lt}.\ref{it:lt4} $G_x$ has
parameter system of length $4$ which is both locked and tied of type
T\ref{it:tied1}.\\
\ee
\begin{figure}
  \begin{center}
    \psfrag{a}{$a$}
    \psfrag{b}{$b$}
    \psfrag{c}{$c$}
    \psfrag{d}{$d$}
    \psfrag{e}{$e$}
    \psfrag{f}{$f$}
    \psfrag{x}{$x$}
    \mbox
    {
      \parbox[b]{.3\textwidth}
      {
        \centerline{
          \includegraphics[scale=0.4]{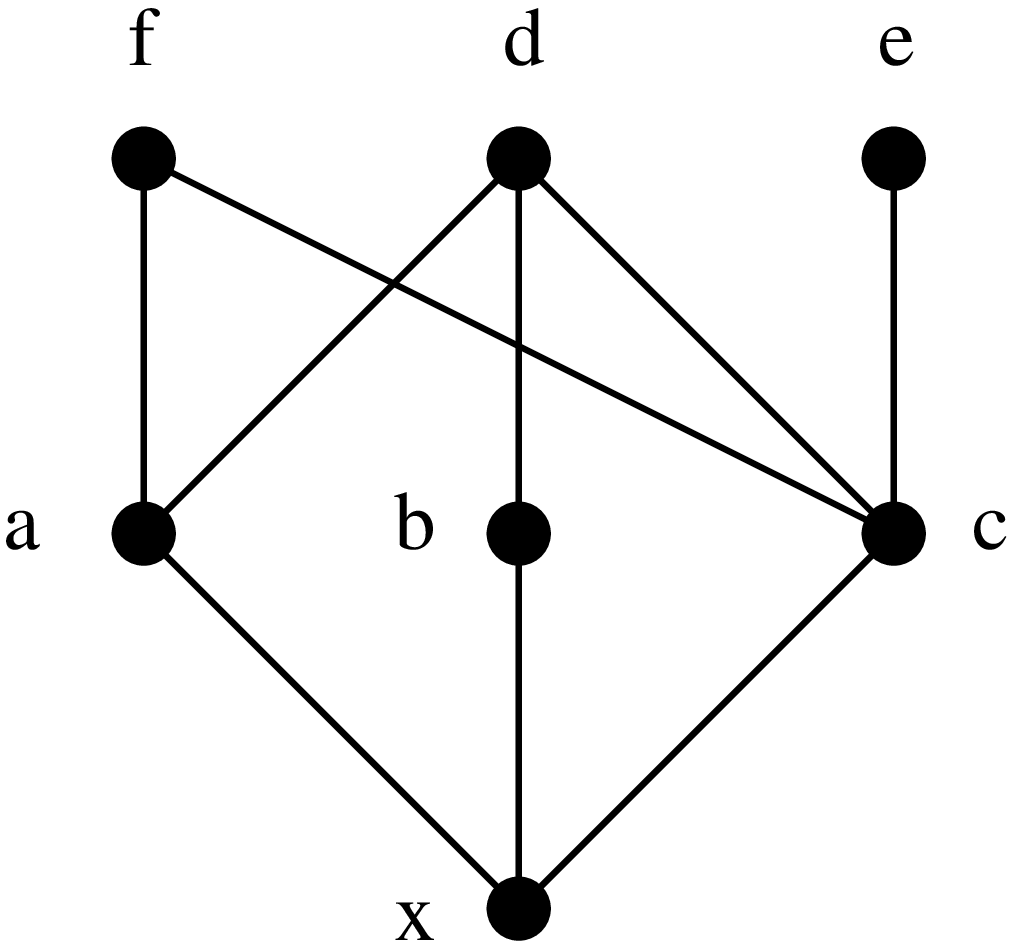} \\[.5em]%
        }
        \centerline{$\G_1$}
      }
      \qquad
      \parbox[b]{.3\textwidth}
      {
        \centerline{
          \includegraphics[scale=0.4]{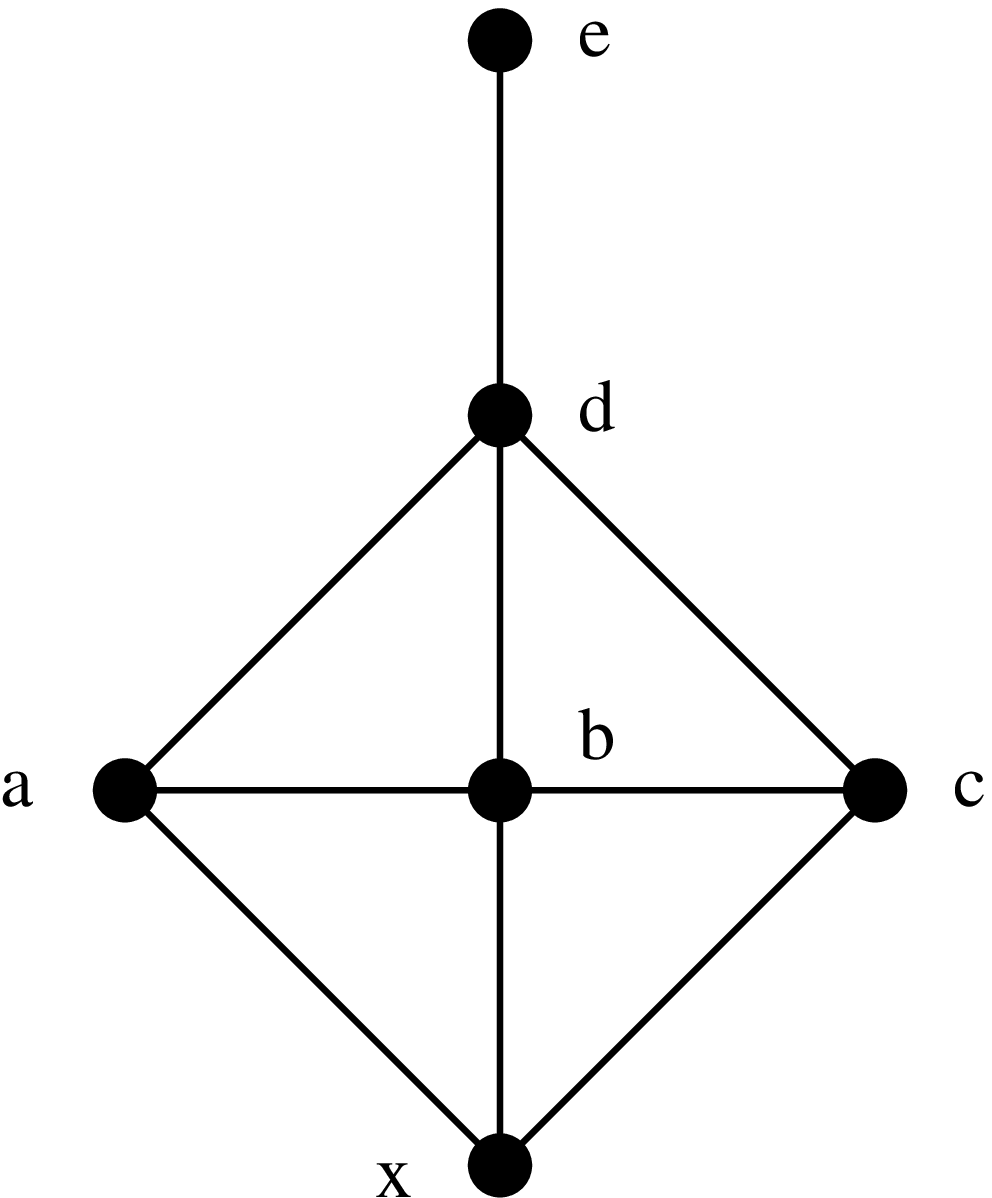}\\[.5em]%
        }
        \centerline{$\G_2$}
      }
      \qquad
      \parbox[b]{.3\textwidth}
      {
        \centerline{\includegraphics[scale=0.4]{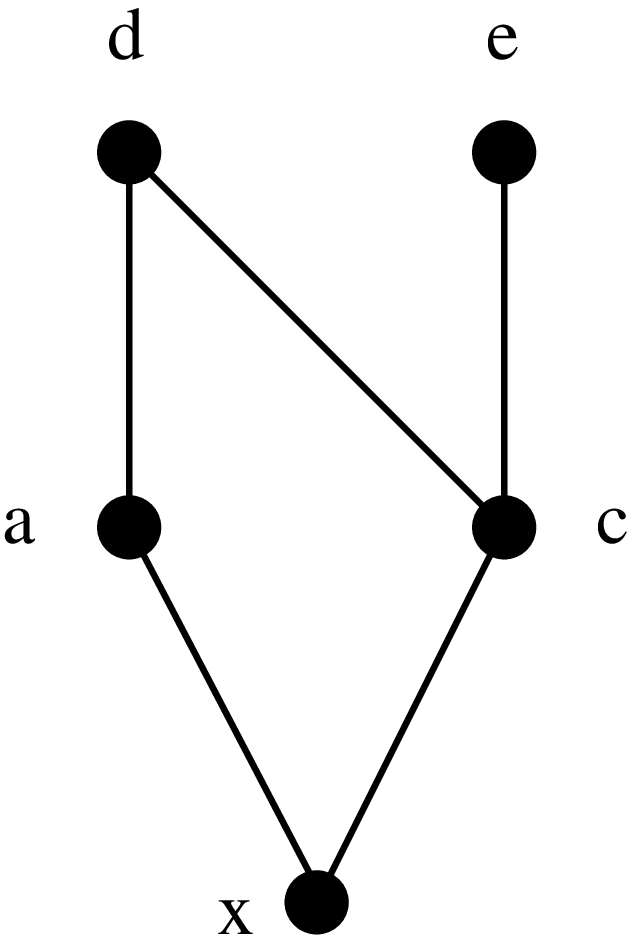}\\[.5em]}
        \centerline{$\G_3$}
      }
    }
  \end{center}
  \caption{Example \ref{ex:lt}}\label{fig:lt1}
\end{figure}
\begin{figure}
  \begin{center}
    \psfrag{a}{$a$}
    \psfrag{b}{$b$}
    \psfrag{c}{$c$}
    \psfrag{d}{$d$}
    \psfrag{e}{$e$}
    \psfrag{f}{$f$}
    \psfrag{g}{$g$}
    \psfrag{h}{$h$}
    \psfrag{i}{$i$}
    \psfrag{j}{$j$}
    \psfrag{k}{$k$}
    \psfrag{l}{$l$}
    \psfrag{x}{$x$}
    \mbox
    {
      \parbox[b]{.3\textwidth}
      {
        \centerline{
          \includegraphics[scale=0.4]{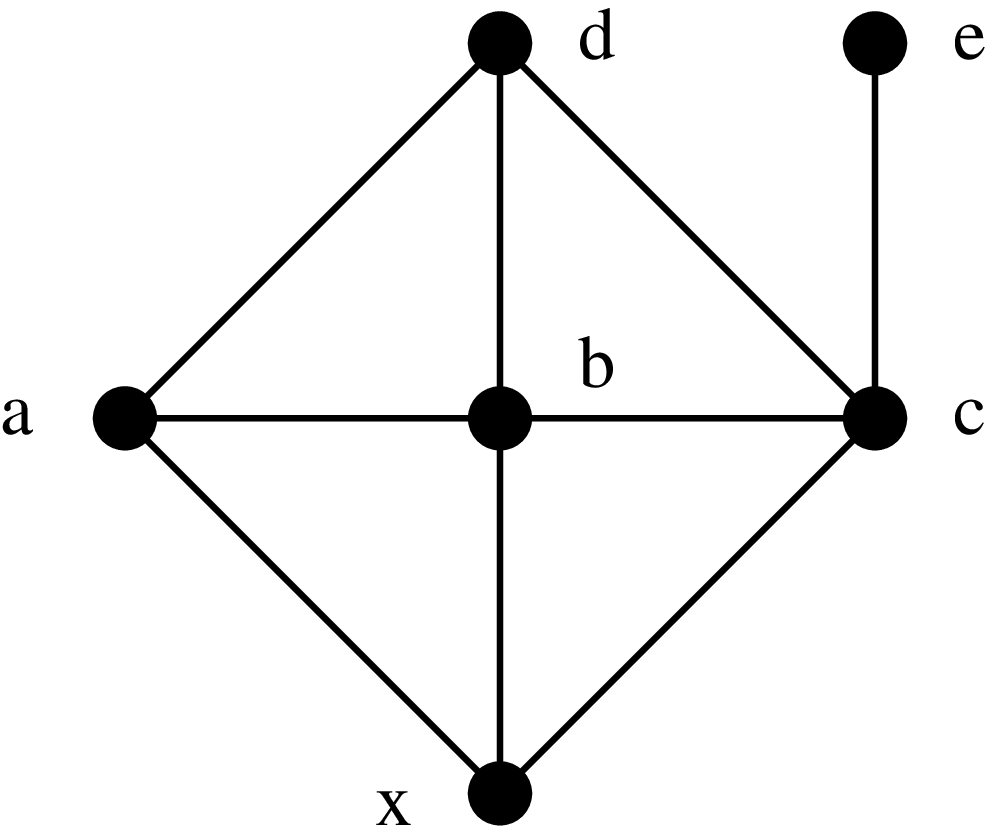} \\[.5em]%
        }
        \centerline{$\G_4$}
      }
      \qquad
      \parbox[b]{.6\textwidth}
      {
        \centerline{
          \includegraphics[scale=0.4]{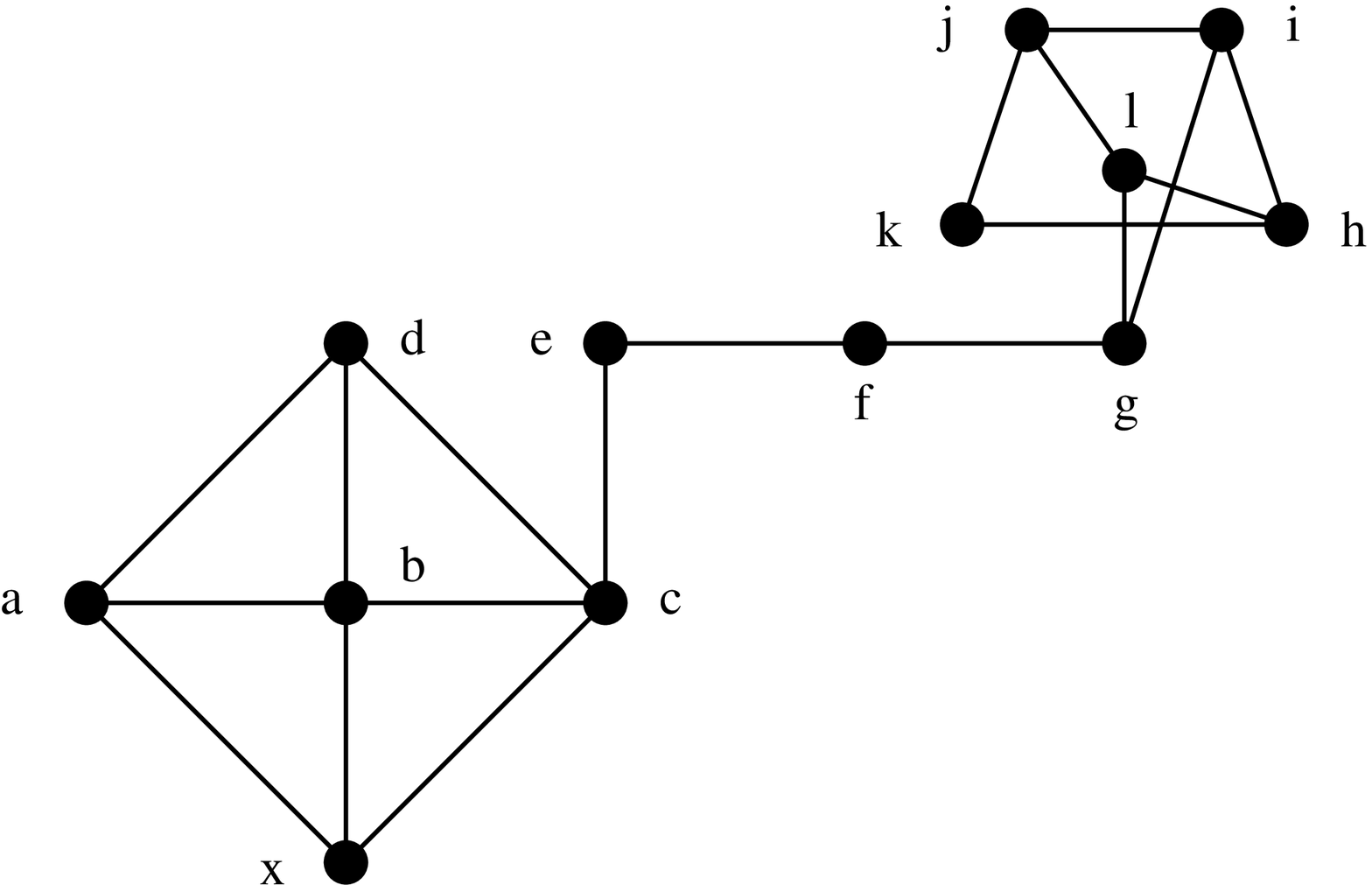}\\[.5em]%
        }
        \centerline{$\G_5$}
      }
    }
  \end{center}
  \caption{Example \ref{ex:lt}}\label{fig:lt2}
\end{figure}
\end{expl}

 Let $P=v_1,\ldots , v_m$ be a maximal parameter system for
$G$.
Let
$D_s=C_G(v_1,\ldots, v_s)$ and let $D_0=G$.
If $x$ occurs in $P$ with  $x=v_j$ then define $j_x=j$. Otherwise set $j_x=m+1$.
 If $x\notin Z(G)$ then define $i_x$ to be the minimal integer $i>0$ such
that $v_i\in W$. Such $i\le m$ exists as $x\notin Z(G)=D_m$;
so there must be some element of $W$ in $P$. If $x\in Z(G)$ then set
$i_x=m+1$.
Note that in this case $W=\nul$ and $x$ does not occur in $P$, so that
$j_x=m+1$. Suppose now that $x$ occurs in $P$ and $j_x=j$. Then
$D_{j-1}\nsubseteq C_G(x)$. In this case suppose that
$D_{j-1}\cap C_G(v_{j+1})\nsubseteq C_G(x)$. As $P$ is a parameter
system
$D_{j-1}\sdc D_{j-1}\cap C_G(v_{j+1})$. There exists $u\in X$ such
that
$u\in (D_{j-1}\cap C_G(v_{j+1}))\backslash C_G(x)$, so $u\neq x$ and we
have $D_{j-1}\cap C_G(v_{j+1})\sdc D_{j-1}\cap C_G(v_{j+1},x)=D_{j+1}$.
Hence we may replace $P$ with the parameter system $v_1, \ldots ,
v_{j-1}, v_{j+1}, x, v_{j+2},\ldots , v_d$. Continuing this way we
eventually replace $P$ with a parameter system such that 
\begin{equation}\label{eq:xmax}
D_{j_x-1}\cap C_G(v_{j_x+1})\subseteq C_G(x).
\end{equation}
Thus, if $x$ occurs in $P$ then $G$ has a maximal parameter system with this
property. Assume now that
$x$ occurs in $P$ and that \eqref{eq:xmax} holds. Let $j_x=j$ and
$i_x=i$.
Suppose that  $i=j+1$ and there is $z\in W$
such that $C_{G_x}(z)\supseteq C_{G_x}(v_1,\ldots , v_{j-1})$. Then we
replace $x$ by $z$ in $P$. This is possible since
$C_{G_x}(z)\supseteq C_{G_x}(v_1,\ldots, v_{j-1})$ implies that
the only element of
$X\cap D_{j-1}$ not in $C_G(z)$ is $x$.
As $D_{j-1}\cap C_G(v_{j+1})\subseteq
C_G(x)$ and $D_{j-1}\nsubseteq C_G(v_{j+1})$
there must be $w\in W$ such that $w\in D_{j-1}$; so $w\in D_{j-1}\cap
C_G(z)$ and $w\notin C_G(v_{j+1})$. Hence after replacing $x$ with $z$
we still have a maximal parameter system for $G$.
In this case we may replace $P$ with a parameter system in which $x$ does
not occur. Thus we may assume that if $P$ contains $x$ and $i_x=j_x+1$ then
there is
no $z\in W$ such that $C_{G_x}(z)\supseteq C_{G_x}(v_1,\ldots,
v_{j_x-1})$. We call a parameter system which
satisfies \eqref{eq:xmax} and the property that either $i_x\neq j_x+1$ or
\[C_{G_x}(z)\nsupseteq C_{G_x}(v_1,\ldots,
v_{j_x-1}), \textrm{ for all } z\in W,\]
a {\em good} parameter system.

Assume now that $P$ is a good maximal parameter system for $G$ with $i_x=i$ and
$j_x=j$.
Then either $x\in Z(G)$ or  $x\in D_{i-1}\backslash C_G(v_i)$.
If either
\be[E1.]
\item\label{it:E1} $x\in Z(G)$ or
\item\label{it:E2} $i=j+1$ or
\item\label{it:E3}  $i\neq j+1$ and $D_{i-1}\backslash C_G(v_i)$ contains $u\in
  X_x$,
\ee
then we say that $P$ {\em extends over} $W$ ({\em at} $i$).
Now
define integers $k$, $l$ and $d$ and a sequence $R=R(P)=y_1,\ldots, y_d$
{\em associated to} $P$ consisting
of
elements
of $X_x$ as follows.
\be[(i)]
\item\label{it:Gxps1} $R=v_1,\ldots ,v_{j-1},v_{j+1},\ldots
  v_{i-1},v_{i+1},\ldots v_d$, $d=m-2$, $k=j-1$ and  $l=i-2$
if $x$ occurs in $P$, $P$ does not
  extend over $W$  and $i>j$.
\item\label{it:Gxps2}   $R=v_1,\ldots ,v_{i-1},v_{i+1},\ldots
  v_{j-1},v_{j+1},\ldots v_d$, $d=m-2$, $k=j-2$ and  $l=i-1$
if $x$ occurs in $P$, $P$ does not
  extend over $W$  and $i<j$.
\item\label{it:Gxps3} $R=v_1,\ldots
  v_{i-1},v_{i+1},\ldots v_d$, $d=m-1$, and  $l=i-1$
if $x$ does not occur in $P$ and  $P$ does not
  extend over $W$.
\item\label{it:Gxps4}   $R=v_1,\ldots
  v_{j-1},v_{j+1},\ldots v_d$, $d=m-1$, and  $k=j-1$
if $x$ occurs in $P$ and  $P$
  extends over $W$.
\item\label{it:Gxps5}   $R=v_1,\ldots
  v_d$, $d=m$,
if $x$ does not occur in $P$ and  $P$
  extends over $W$.
\ee
\begin{lem}\label{lem:R(P)}
Let $P$ be a good maximal parameter system for $G$.
Then $R=R(P)$ is a parameter system of length $d$ for $G_x$, where
$0\le \cd(G) -d\le 2$. Moreover exactly one of the following conditions
R1--R4 holds.
\be[R1.]
\item\label{it:R1} $R$ is locked and tied, of type T\ref{it:tied1} but not
T\ref{it:tied2}, $d=m-2$, $x$
  occurs in $P$ and $P$ does not extend over $W$.
\item\label{it:R2} $R$ is locked $d=m-1$, $x$ does not
  occur in $P$ and $P$ does not extend over $W$.
\item\label{it:R3} $R$ is tied, $d=m-1$, $x$
  occurs in $P$ and $P$ extends over $W$.
\item\label{it:R4} $d=m$,  $x$ does not
  occur in $P$ and $P$ extends over $W$.
\ee
\end{lem}
\begin{proof}
First we show that $R$ is a parameter system for $G_x$.
If $1\le s\le m$ and $s\neq i$ then there exists $u\in X_x$ such
that $u\in D_{s-1}\backslash D_s$. Therefore in cases \eqref{it:Gxps1} to
\eqref{it:Gxps3} $R$ is a parameter system for $G_x$.
If $P$ extends over $W$ at $i$, $x\notin Z(G)$  and $i\neq j+1$ then
there exists $u\in X_x$ such
that $u\in D_{i-1}\backslash D_i$.
Hence the required result holds in case \eqref{it:Gxps4} when
$i\neq j+1$.
Suppose that $P$ extends over $W$, $x\notin Z(G)$  and $i=j+1$.
Then $x=v_j=v_{i-1}$
 so $C_G(v_1,\ldots,v_{i-2})\nsubseteq C_G(x)$ and
$C_G(v_1,\ldots,v_{i-2})\cap C_G(v_i)\subseteq C_G(x)$.
Hence there exists
$w\in X\backslash C_G(x)$ such that $w\in C_G(v_1,\ldots ,v_{i-2})$,
with $w\notin C_G(v_i)$. Thus $w\in X_x$ and
$C_{G_x}(v_1,\ldots ,v_{i-2})\sdc C_{G_x}(v_1,\ldots,v_{i-2})\cap
C_{G_x}(v_i)$.
It follows that $R$ is a parameter
system for $G_x$ in case \eqref{it:Gxps4}.
In case \eqref{it:Gxps5} if $x\in Z(G)$ then $W=\nul$ so
$P=R$ is a parameter system for $G_x$. Otherwise, in case
\eqref{it:Gxps5}, $i\neq j+1$ and
the same arguments
as in case \eqref{it:Gxps4} apply.

Clearly $0\le \cd(G)-d\le 2$ so it remains to prove that exactly one
of R\ref{it:R1} to R\ref{it:R4} holds.
If $P$ does not extend
over $W$ then $i\neq j+1$ and
\begin{equation}\label{eq:notext}
D_{i-1}\backslash C_G(v_i)\cap X=\{x\}.
\end{equation}
In this case suppose first that $i<j$ (so we are in case \eqref{it:Gxps2}
or \eqref{it:Gxps3}).
Then $v_i\in W$, $l=i-1$ and $x\notin \{ v_1,\ldots ,v_{i-1}\}\subseteq Y$
and \eqref{eq:notext} implies
$C_{G_x}(y_1,\ldots,y_l)\subseteq C_{G_x}(v_i)$. Thus
$R$ is locked at $l$ with key $v_i$.

Now suppose that
$P$ does not extend over $W$ and that $i>j$ (so we have case
\eqref{it:Gxps1}). This time
$D_{i-1}=C_G(v_1,\ldots, x,\ldots,v_{i-1})$ and
\eqref{eq:xmax} implies that
$C_G(v_1,\ldots ,v_{j-1},v_{j+1})\subseteq C_G(x)$.
Thus \eqref{eq:notext} reads
\[
[C_G(v_1,\ldots ,v_{j-1},v_{j+1},\ldots ,v_{i-1})\cap C_G(x)]
\backslash C_{G}(v_i)\cap
X=\{x\}\]
so that
\[
[C_G(v_1,\ldots ,v_{j-1},v_{j+1},\ldots ,v_{i-1})]
\backslash C_{G}(v_i)\cap X=\{x\}
\]
which implies
$C_{G_x}(y_1,\ldots, y_l)\subseteq G_x(v_i)$.
Again $R$ is locked at $l$ with key $v_i$.
Therefore, if $P$ does not extend over $W$ we have a locked parameter
system $R$ for $G_x$. The length of $R$ is $d = \cd(G)-2$, if $x$
occurs in $P$, and
$d = \cd(G)-1$ otherwise.
On the other hand, considering cases \eqref{it:Gxps1}-\eqref{it:Gxps5},
if $P$ extends over $W$ then $R$ has length
$d\ge \cd(G)-1$.

Now suppose that $x$ occurs in $P$.
Then $x=v_j$  and
$D_{j-1}\nsubseteq C_G(x)$ so
\begin{equation}\label{eq:knotA}
C_{G_x}(y_1,\ldots, y_k)\nsubseteq A.
\end{equation}
Firstly suppose that $i\neq j+1$ and $j\neq m$.
We claim that in this case
\begin{equation}\label{eq:inA}
C_{G_x}(y_1,\ldots ,y_k)\cap C_{G_x}(y_{k+1})\subseteq A.
\end{equation}
To see this note that if
$i>j$ or $P$ extends over $W$ then $\{y_1,\ldots ,y_k\}=
\{v_1,\ldots ,v_{j-1}\}$ (see cases \eqref{it:Gxps1} and \eqref{it:Gxps4}).
Therefore in these cases, as $y_{k+1}=v_{j+1}$, \eqref{eq:xmax} implies that
\eqref{eq:inA} holds. If $P$ does not extend over $W$ and $i<j$ then $k>l$ and
$\{y_1,\ldots ,y_l\}=\{v_1,\ldots ,v_{i-1}\}$. Thus if
$u\in X_x\cap C_{G_x}(y_1,\ldots ,y_k)$ and $u\notin C_G(v_i)$ then
$u\in D_{i-1}\backslash C_G(v_i)$, so $P$ extends over $W$, a contradiction.
Hence such $u$ must belong to $C_G(v_i)$ and therefore to $D_{j-1}$. Now
\eqref{eq:xmax} implies that if in addition $u\in C_{G_x}(y_{k+1})$
then $u\in C_G(x)$ so $u\in A$. Hence \eqref{eq:inA} holds
in this case as well.
As $D_j\subseteq C_G(x)$ and $D_j\sdc
D_{j+1}$ there exists $y\in Y\cup \{x\}$ such that $y\in D_j$ and
$y\notin D_{j+1}$. As $i=j+1$
we have $y\neq x$, so $y\in Y$ and therefore $y\in C_{G_x}(y_1,\ldots
,y_k)\backslash C_{G_x}(y_1,\ldots, y_{k+1})$. As $i\neq j+1$ and $j\neq m$
we have $k<n$ and $y_{k+1}\notin W$. Hence $R$
is tied at $k$  of type T\ref{it:t1a} (and not of type T\ref{it:t1b} or
T\ref{it:tied2}).

Secondly suppose that $j\neq m$, $i=j+1$. Then $P$ extends over $i$,
$k=j-1$ and $y_1,\ldots,y_k=v_1,\ldots, v_{j-1}\in Y$
(case \eqref{it:Gxps4}). As before \eqref{eq:knotA} holds and \eqref{eq:xmax}
implies \eqref{eq:inA}.
This time $y_{k+1}=v_{j+1}\in W$ so
$R$ is tied at $k$
of type T\ref{it:tied2} (and possibly also of type T\ref{it:t1a}).
Since $P$ is good
$R$ is not locked.

Thirdly suppose that $j=m$; in which case $i\neq j+1$. Then $d=k$ and
there exists $z\in X$ such that $z\in C_G(v_1,\ldots ,v_{m-1})$ but
$z\notin C_G(x)$. Hence $z\in W$ and
$z\in C_{G_x}(y_1,\ldots,y_d)$.
Therefore $R$ is tied at $k$ of type T\ref{it:t1b} (but not T\ref{it:tied2}).

Finally note that R\ref{it:R4} occurs if and only if
we have case \eqref{it:Gxps5}. If $P$ does not extend over
$W$ and $x$ occurs in $P$ then, from the above, $R$ is locked at
$l$ and tied at $k$ of type T\ref{it:t1a} or T\ref{it:t1b}, but not T\ref{it:tied2}.
In this case R\ref{it:R1} is
satisfied. If $P$ does not extend over $W$ and $x$ does not occur in $P$ then
$R$ is locked at $l$ and R\ref{it:R2} holds. If $P$ extends over $W$ and $x$ occurs in
$P$ then R\ref{it:R3} holds.
\end{proof}

Let $Q=x_1,\ldots,x_n$ be a parameter system for $G_x$. Let
$C_s=C_{G_x}(x_1,\ldots ,x_s)$,
and
$C_0=G_x$.
Suppose that $Q$
is locked at $l$. Note that if $l<l^\prime\le n$ and $x_i\in Y$, for
$i=1,\ldots ,l^\prime$, then $Q$ is also locked at $l^\prime$.
If $Q$ is locked at $l$ and $l$ is maximal such that $x_i\in Y$, for
$i\le l$, then
we say that $Q$ is {\em right} locked at $l$.
We may always assume that if $Q$ is locked at $l$
then it  is right locked at $l$.
Note also that if $Q$ is tied of type T\ref{it:tied2} and $Q$ is
right locked  at $l$ then $x_{l+1}\in W$
so that $k=l<n$. Conversely, if $Q$ is both right locked and tied at $l<n$
then
it is tied of type T\ref{it:tied2} (and possibly also of type
T\ref{it:tied1}).

Now suppose that $Q$ is either right locked at $l$ with key $w$ or tied at
$k$, or both.  Define  integers $d$, $i$ and $j$  and a
sequence $S=S(Q)=u_1,\ldots,u_d$
of elements of $G$ as follows.
\be[(i)]
\item $S=x_1,\ldots ,x_k,x,x_{k+1},\ldots ,x_l,w,x_{l+1},\ldots ,x_n$,
  $d=n+2$, $i=l+2$ and $j=k+1$, if either
  \be[I.] 
\item $Q$ is locked and tied of type
  T\ref{it:tied1} and $k<l$; or 
\item $Q$ is locked and tied of type  T\ref{it:t1b} and $k=l=n$
  .
  \ee
\item $S=x_1,\ldots ,x_l,w,x_{l+1},\ldots ,x_k,x,x_{k+1},\ldots ,x_n$,
  $d=n+2$, $i=l+1$ and $j=k+2$, if $Q$ is locked and tied of type
  T\ref{it:tied1} and $k>l$.
\item\label{it:S3} $S=x_1,\ldots ,x_l,w,x_{l+1},\ldots ,x_n$,
  $d=n+1$, $i=l+1$, if $Q$ is locked but not tied.
\item\label{it:S4}  $S=x_1,\ldots ,x_k,x,x_{k+1},\ldots ,x_n$,
  $d=n+1$, $j=k+1$, if either $Q$ is not locked but is tied of type T\ref{it:tied1} or $Q$ is
  tied of type T\ref{it:tied2}.
\item\label{it:S5} $S=x_1,\ldots ,x_n$, $d=n$,
  if $Q$ is neither locked nor  tied.
\ee
\begin{lem}\label{lem:S(Q)}
Let $Q$ be a parameter system for $G_x$. Then $S(Q)$ is a parameter
system for $G$.
\end{lem}
\begin{proof}
Let $E_s=C_G(u_1,\ldots u_s)$, for $s=1,\ldots ,d$, and $E_0=G$.
First suppose that $Q$ is locked and tied of type
  T\ref{it:tied1} and $k<l$. Clearly
$E_0\sdc E_1\sdc \cdots \sdc E_{j-1}$. As
  $Q$ is tied at $k$ with $k<l\le n$ it must be tied of type
  T\ref{it:t1a}, so $k$ is maximal such that $C_k\nsubseteq A$ and
there exists $y\in Y$ and $z\in W$ such that $y,z\in C_k$ but
$y,z\notin C_{k+1}$. Hence $y,z\in E_{j-1}$, $y\in E_j$, $z\notin E_j$
and $y\notin E_{j+1}$; so $E_{j-1}\sdc E_j\sdc E_{j+1}$. Since
$C_{k+1} \subseteq A$ we have $C_G(x_1,\ldots,x_k)\cap C_G(x)\cap
C_G(x_{k+1},\ldots ,x_{t})= C_G(x_1,\ldots,x_k)\cap
C_G(x_{k+1},\ldots ,x_{t})$, for $t=k+1,\ldots, i-1$, so we also have
$E_{j+1}\sdc \cdots \sdc E_{i-1}$.
As $Q$ is locked at $l$ with key $w$ we have $x_i\in Y$, for
$i=1,\ldots, k,k+1,\ldots, l$, and $C_{G_x}(w)\supseteq C_l$.
Hence $x\in C_G(u_1,\ldots, u_{i-1})$ and $x\notin  C_G(u_1,\ldots,
u_{i-1},w)$. Moreover, if $z\in X$ and $z\in  C_G(u_1,\ldots,
u_{i-1})\backslash C_G(w)$ then $z=x$.
As $Q$ is right locked at $l$,
$u_{i+1}=x_{l+1}\in W$ so
$C_G(u_1,\ldots,
u_{i-1},w,u_{i+1},\ldots ,u_t)=C_G(u_1,\ldots,
u_{i-1},u_{i+1},\ldots ,u_t)$, for all $t\ge i+1$. It follows that
$E_{i-1}\sdc E_i\sdc E_{i+1}\sdc \cdots \sdc E_d$ and so $S$ is a
parameter system for $G$.

If $Q$ is locked and tied of type
  T\ref{it:t1b} and $k=l=n$ then $E_0\sdc \cdots \sdc E_n$. As $x,w\in
  E_n$, $x\in E_{n+1}\backslash E_{n+2}$ and $w\notin E_{n+1}$ it
  follows
that $S(Q)$ is a parameter system for $G$.

Now suppose that $Q$ is locked and tied of type
  T\ref{it:tied1} and $k>l$. As above
$E_1\sdc \cdots \sdc E_{i-1}\sdc E_i\sdc E_{i+1}\sdc \cdots \sdc
E_{j-1}$ and $E_{j-1}=C_G(x_1,\ldots ,x_l)\cap C_G(x_{l+1},\ldots ,x_k)$.
First assume that $k\neq n$, so that $Q$ is tied of type
T\ref{it:t1a}.
By definition of tied there is
$y\in Y$ and $z\in W$ such that $\{y,z\}\subseteq C_k\backslash
C_{k+1}$ and, from the above, $\{y,z\}\subseteq E_{j-1}$.
As in the previous case it now follows that $S$
is a parameter system for $G$. In the case where $k=n$ we have
$E_1\sdc \cdots \sdc E_{d-1}$ and there exists $z\in W$ such that
$z\in C_n=Z(G_x)$. Hence 
$z\in E_{d-1}$ but
$z\notin E_d$. Again $S$  is a parameter system for $G$.

Similar
arguments show that in case (\ref{it:S3}) and, in the instances of
case (\ref{it:S4})
when $Q$ is tied of type T\ref{it:tied1}, $S$ is a parameter system
for $G$. Obviously $S$ is a parameter system for $G$ in case \eqref{it:S5}.
This leaves the case where $Q$ is tied of type T\ref{it:tied2} at $k$.
Note that if also $Q$ is right locked at $l$ then $k=l$.
Then $E_1\sdc \cdots \sdc E_{j-1}$ and there exists $z\in W$ such that
$z\in C_k$ but $z\notin C_{k+1}$. Since $C_{k+1}\subseteq A$ we have
$C_G(x_1,\ldots ,x_{k+1})\subseteq C_G(x)$ and $z\in W$ implies
$z\notin C_G(x)$. As $x_s\in Y$, for $s=1,\ldots k$ we have $x\in
C_G(x_1,\ldots ,x_k)$. Hence $z,x\in E_{j-1}$, $x\in E_j$,
$z\notin E_j$ and
$x\notin E_{j+1}$, as
$x_{k+1}\in W$,  so $E_{j-1}\sdc E_j\sdc E_{j+1}$. As before
$E_{j+1}\sdc \cdots \sdc E_d$, so $S$ is a parameter system for $G$.
Thus in all cases $S$ is a parameter system for $G$.
\end{proof}
If a parameter system for $G_x$ is tied of type T\ref{it:tied1} but is not tied
of type T\ref{it:tied2} then we say that it is of
type T\ref{it:tied1}-{\em exclusive}.
\begin{thm}\label{cd-ext}
  In the notation above $0\le\cd(G)-\cd(G_x)\le 2$ and, more precisely, the
  following hold.
  \be[\rm 1.]
\item $\cd(G)-\cd(G_x)=2$ if and only if $G_x$ has a maximal
  parameter system which is both locked and tied of type
  T\ref{it:tied1}-exclusive.
\item $\cd(G)-\cd(G_x)=1$ if and only if $G_x$ has a parameter
  system which is either
  \be
\item maximal and locked or tied but not both locked and tied
of type T\ref{it:tied1}-exclusive; or
\item of length $\cd(G_x)-1$, ends in $Z(G_x)$,  and is
both locked and tied of type T\ref{it:tied1}-exclusive.
  \ee
\item Otherwise $\cd(G)=\cd(G_x)$.
  \ee
\end{thm}
\begin{proof}
As $G_x$ is a subgroup of $G$ it follows that $\cd(G_x)\le
\cd(G)$.
Suppose that $G$ has a maximal parameter system $Q$ which is locked at
$l$
and tied, of type T\ref{it:tied1}-exclusive, at $k$.
As above, $l\neq k$. From Lemma \ref{lem:S(Q)} $S(Q)$
is a parameter system for $G$ so
$\cd(G)\ge \cd(G_x)+2$ and, from Lemma \ref{lem:R(P)},
we conclude that $\cd(G)= \cd(G_x)+2$.

Conversely, suppose that $\cd(G)= \cd(G_x)+2$ and let $P$ be a maximal
parameter system for $G$. We may assume that $P$ is good. Then if
$R=R(P)$ satisfies R\ref{it:R2}, R\ref{it:R3} or R\ref{it:R4} we have
$\cd(G)\le
\cd(G_x) +1$, a contradiction. Hence $R$ satisfies R\ref{it:R1} and
since $R$
has length $\cd(G)-2=\cd(G_x)$ it must be maximal. This
proves the first numbered statement of the Theorem.

Now suppose that $\cd(G)\le \cd(G_x)+1$; so, from Lemma \ref{lem:S(Q)}
$G_x$ has no maximal
parameter system which is both locked and tied of type
T\ref{it:tied1}-exclusive.
Suppose that $G_x$ has a maximal parameter system $Q$ which is either locked
or tied. From Lemma \ref{lem:S(Q)}, $S(Q)$ is a parameter system
for $G$ of
length
$\cd(G_x)+1$. Hence, in this case $\cd(G)=\cd(G_x)+1$. If $G_x$ has a
parameter
system $Q$ of length $\cd(G_x)-1$ which is both locked and tied of
type T\ref{it:tied1}-exclusive then $S(Q)$ has length $\cd(G_x)+1$,
so again
$\cd(G)=\cd(G_x)+1$.

Conversely suppose that $\cd(G)=\cd(G_x)+1$, let $P$ be a good
maximal parameter system for $G$ and let $R=R(P)$.
Then $R$ does not satisfy R\ref{it:R4},
since this  would imply $\cd(G_x)\ge \cd(G)$. From Lemma
\ref{lem:R(P)}, $R$ satisfies one of R\ref{it:R1},
R\ref{it:R2} or R\ref{it:R3}. If $R$ satisfies
R\ref{it:R2} or
R\ref{it:R3}
then $R$ has length $\cd(G)-1=\cd(G_x)$ so is maximal. In this
case
$G_x$ has a maximal parameter system which is either locked or tied.
If $R$ satisfies R\ref{it:R1} then $R$ has length $\cd(G)-2=\cd(G_x)-1$.
If $R$ does not end in $Z(G_x)$ then there exists $u\in X_x$ such that
$R^\prime=y_1,\ldots,y_d,u$ is a parameter system for $G_x$, which must be maximal.
Since $R$ is locked at $l\le d$, $R^\prime$ is also locked at $l$ and so $G_x$
has a maximal parameter system which is locked. Thus we may assume that $R$
ends in $Z(G_x)$.
This
proves the second numbered statement of the theorem and the third follows
from Lemma \ref{lem:R(P)}.
\end{proof}

In some cases we have simple criteria which determine the relationship
of $\cd(G)$ to $\cd(G_x)$. We sum these up in the next two
corollaries. Observe though that, of the extensions of Example
\ref{ex:lt}, only \ref{it:lt2} falls within the scope of Corollary
\ref{ext-simple}.

\begin{cor}\label{ext-simple} In the notation above let $B=G(W)$.
\begin{enumerate}
\item\label{ext-simple1} $\cd(G)=\cd(G_x)+2$ if
\begin{enumerate}
\item\label{ext-simple1a} $Z(G)\sac Z(G_x)$; or
\item\label{ext-simple1b} $\cd(A)=\cd(G_x)$.
\end{enumerate}
\item\label{ext-simple2} $\cd(G)=\cd(G_x)$ if
\begin{enumerate}
\item\label{ext-simple2a} $G_x=A\times B$ and $Z(B)=1$; or
\item\label{ext-simple2b} $C_G(x)=C_G(S)$, for some $S\subseteq G_x$.
\end{enumerate}
\end{enumerate}
\end{cor}
\begin{proof}  Let $Q=x_1,\ldots ,x_n$ be parameter system for $G_x$ and
  let $C_s=C_{G_x}(x_1,\ldots, x_s)$.
\be
\item
\be
\item If $Z(G)\sac Z(G_x)$ then there exists $z\in W$ such that $z\in
  Z(G_x)$.
Suppose that $Q$ is a
  maximal parameter system for $G_x$. Then $Q$ is tied of type T\ref{it:t1b} at $n$,
so is tied of type T\ref{it:tied1}-exclusive.  Let $l$ be
maximal such that $x_s\in Y$, for all $s\le l$. Since $z\in Z(G_x)$ we
have $C_{G_x}(z)\supseteq C_l$. Hence $Q$ is locked at $l$ and the
result follows from Theorem \ref{cd-ext}.
\item We may now assume that $Z(G)=Z(G_x)$. Removing all generators in
  the centre of $G$ results in a new graph group with the same
  centraliser
dimension as $G$. Therefore we may assume that $Z(G)=Z(G_x)=1$. Let
$T=t_1,\ldots ,t_d$ be a maximal parameter system for $A$. This is
also
a parameter system for $G_x$, which must be maximal since
$\cd(A)=\cd(G_x)$. Let
$F_s=C_{G_x}(t_1,\ldots ,t_s)$, for all $s$. Then $F_d=Z(G_x)=1$ and so there
exists $w\in W$ such that $C_{G_x}(w)\supseteq F_d$. As $t_i\in A$ we
have $t_i\in Y$, for all $i$, so as a parameter system for $G_x$, $T$ is locked at
$d$ with key $w$. As $T$ is a parameter system for $A$ there is
$y_s\in Y$ such that $y_s\in F_{s-1}\backslash F_s$, for $s=1,\ldots
d$. Let $k$ be maximal such that $F_k\nsubseteq A$. As $Z(G_x)=1$, we
have $k<d$. Then $T$ is tied
of type T\ref{it:t1a} at $k$. From Theorem \ref{cd-ext} then
$\cd(G)=\cd(G_x)+2$.
\ee
\item
\be
\item Suppose that $Q$ is a parameter system for $G_x$ which ends in $Z(G_x)$.
Let $l$ be maximal such that $x_s\in Y$, for $s=1,\ldots
  ,l$. Then $W\subseteq C_l$ and so if there is $w\in W$ such that
  $C_{G_x}(w)\supseteq C_l$ then $C_{G_x}(w)\supseteq W$,
  contradicting the hypothesis that $Z(B)=1$. Hence $Q$ is not
  locked. Let
$k$ be maximal such that $C_k\nsubseteq A$. Then, as $Z(G_x)=1$ we
have
$k<n$. Every element of $W$ commutes with every element of $Y$, so
$C_{k+1}\subseteq A$ implies $x_{k+1}\in W$. Therefore
$C_{G_x}(x_{k+1})\supseteq Y$ and $Q$ cannot be tied of type
T\ref{it:tied1}. If $Q$ is tied of type T\ref{it:tied2} then
$x_s\in Y$, for $s=1,\ldots k$, in which case $B\subseteq
C_k$. Therefore $x_{k+1}\in C_{k+1}$, a contradiction. Hence
$G_x$ has no parameter system which is either locked or
tied and ends in $Z(G)$. From Theorem \ref{cd-ext}, $\cd(G)=\cd(G_x)$.
\item Since $C_G(x)=C_G(S)$ we have $A=C_{G_x}(S)$, so that $y\in Y$
  implies $S\subseteq C_{G_x}(y)$. If $x_t\in Y$, for $t=1,\ldots ,l$,
  then $S\subseteq C_l$. Moreover, if $w\in W$ then
  $[w,s]\neq 1$, for some $s\in S$, so $S\nsubseteq C_{G_x}(w)$.
  Therefore $Q$ is not locked. Let $k$ be maximal such that
  $C_k\nsubseteq A$. Then $C_{k+1}\subseteq A=C_{G_x}(S)$.
Thus there is at least one $w\in W$ such that $w\in C_k\backslash
C_{k+1}$. Fix such a $w$. Then $w\notin A$ so there exists some $u\in
X_x$ such that $u \in \a(S)$ and $[u,w]\neq 1$.
As $C_{G_x}(u)\supseteq C_{G_x}(S)$ we have
$C_k>C_k\cap C_{G_x}(u)\supseteq C_k\cap C_{G_x}(S)\sdc C_{k+1}$. Hence we may
assume that $x_{k+1}=u$. Now if $y\in Y$ then $y\in C_{G_x}(S)$
implies $y\in C_{G_x}(u)$ so that  $y\notin C_k\backslash
C_{k+1}$. Hence $Q$ is not tied of type T\ref{it:t1a}. As $Z(G)=1$, $Q$
is not tied of type T\ref{it:t1b}. As $u$ occurs as a letter in an
element
of $S$ and $C_G(x)=C_G(S)$ it follows that $u\in A$, so $Q$
is not tied of type T\ref{it:tied2}.  From Theorem \ref{cd-ext}, $\cd(G)=\cd(G_x)$.
\ee
\ee
\end{proof}

\begin{expl}
Using Theorem \ref{cd-ext} we may compute the centraliser dimension
of the semibraid groups. The $n$th {\em semibraid} group is defined to
be the partially commutative
group $G_n$ with
presentation
\[\la x_1,\ldots ,x_n \,|\, [x_i,x_j]=1, \textrm{ for all } i, j
\textrm{ such that } |i-j|\ge 2\ra .
\]
By inspection, using Proposition \ref{thm:pcdimgen},
 it is found that the centraliser dimensions of the first
few semibraid groups are $\cd(G_1)=0$, $\cd(G_2)=\cd(G_3)=2$ and $\cd(G_4)=\cd(G_5)=4$.
Proposition \ref{thm:pcdimgen} also implies that $Z(G_n)=1$, for all $n\ge 2$, as
$C(x_1,\ldots ,x_n)=1$. The following lemma shows that this pattern continues.

\begin{lem}
For $n\ge 0$, $\cd(G_{2n})=\cd(G_{2n+1})=2n$.
\end{lem}
\begin{proof}
Let $n\ge 4$ and suppose that $\cd(G_n)=d$. Let $\G_{n+2}$ be the 
commutation graph of the group $G_{n+2}$. Let $K$ be the partially group
with commutation graph 
obtained from $\G_{n+2}$ by deleting the vertex $x_{n+1}$. Then
$[x_j,x_{n+2}]=1$, for $j=1,\ldots n$, so it follows from Corollary
\ref{ext-simple}.\ref{ext-simple2b}, with $S=1$, that
$\cd(K)=\cd(G_n)$. As $Z(G_{n+2})=1$ and $Z(K)=\la x_{n+2}\ra$ 
Corollary \ref{ext-simple}.\ref{ext-simple1a} implies that
$\cd(G_{n+2})=\cd(K)+2=d+2$.
Thus the lemma follows from the facts that $\cd(G_4)=\cd(G_5)=4$.
\end{proof}
\end{expl}
The example above shows that there are partially commutative groups
of any given even centraliser dimension. 
It is well known that there are no groups of centraliser dimension $1$
(see for example \cite{DKR}). It was pointed out to us by K. Goda
that
there are no partially commutative groups of centraliser dimension $3$
either. In fact if we suppose that such groups do exist then we
can take $G$ to be  a partially commutative 
group of centraliser dimension $3$ with the minimal number of vertices
among all such groups. The centre of $G$ will then be trivial (as
otherwise
we may remove all central generators). Now remove a vertex $x$ from
the commutation graph of $G$ and let $G_x$ be the corresponding group.
Then $\cd(G_x)\neq 1$ and $\cd(G_x)\neq 3$, by minimality of $G$, so
$\cd(G_x)=2$. If $Z(G_x)\neq 1$ then Corollary
\ref{ext-simple}.\ref{ext-simple1a} 
implies that 
$\cd(G)=4$. Hence $Z(G_x)=1$ and so, from \cite[Proposition 3.9.1]{DKR},
$G_x$ is a CT-group. It follows that $G_x$ is a free product of free
Abelian groups. Using Corollary \ref{ext-simple} again it is not hard
to
see that the
centraliser
dimension of $G$ must either be $2$ or $4$, a contradiction.
Therefore there is no partially commutative group of centraliser
dimension $3$.
However, as the following example shows, every odd number greater
than $3$ is the centraliser dimension of a partially commutative group.

\begin{expl}
If $G$ is a partially commutative group of centraliser dimension 
$d$ then it follows from \cite[Lemma 2.2.2]{MS} that the partially
commutative group $G\times F$, where $F$ is free of rank $2$, has
centraliser
dimension $d+2$.
The group of Example \ref{ex:lt} with commutation graph $\G_1$ has
centraliser dimension $5$. Therefore there exist partially commutative
groups of centraliser dimension $2m+1$, for all $m\ge 2$.
\end{expl}

\begin{flushleft}
  {
    \small Andrew J. Duncan,\\
    School of Mathematics and Statistics,\\
    Merz Court, University of Newcastle,\\
    Newcastle upon Tyne, NE1 7RU, United Kingdom.\\
    e-mail: A.Duncan@ncl.ac.uk
  }
\end{flushleft}
\begin{flushleft}
  {
    \small 
    Ilya V. Kazachkov\newline
  Department of Mathematics and Statistics\\
University of McGill\\     
805 Sherbrooke West, Montreal,\\
     Quebec, Canada, H3A 2K6.\\
  e-mail: kazachkov@math.mcgill.ca   
  }
\end{flushleft}
\begin{flushleft}
  {
    \small 
       Vladimir N. Remeslennikov\newline
    e-mail: remesl@iitam.omsk.net.ru\\[.5em]
    Institute of Mathematics \newline 
    (Siberian Branch of Russian Academy of
    Science)\newline  
    Pevtsova st. 13\newline
    644099, Omsk, Russia, \newline 
    +7 381 2 247041 Fax +7 381 2 234584
  }
\end{flushleft}
\end{document}